# Global bifurcation of solitary waves for the Whitham equation

Tien Truong[1] · Erik Wahlén[1] · Miles H. Wheeler[2]

## Abstract

The Whitham equation is a nonlocal shallow water-wave model which combines the quadratic nonlinearity of the KdV equation with the linear dispersion of the full water wave problem. Whitham conjectured the existence of a highest, cusped, traveling-wave solution, and his conjecture was recently verified in the periodic case by Ehrnström and Wahlén. In the present paper we prove it for solitary waves. Like in the periodic case, the proof is based on global bifurcation theory but with several new challenges. In particular, the small-amplitude limit is singular and cannot be handled using regular bifurcation theory. Instead we use an approach based on a nonlocal version of the center manifold theorem. In the large-amplitude theory a new challenge is a possible loss of compactness, which we rule out using qualitative properties of the equation. The highest wave is found as a limit point of the global bifurcation curve.

## 1 Introduction

In this paper, we continue the story of singular wave phenomena featured in the Whitham equation. The equation was proposed by Whitham [33], in an attempt to remedy the failure of the KdV equation in capturing wave breaking and peaking. He

✉ Tien Truong
  tien.truong@math.lu.se

  Erik Wahlén
  erik.wahlen@math.lu.se

  Miles H. Wheeler
  mw2319@bath.ac.uk

[1] Centre for Mathematical Sciences, Lund University, Lund, Sweden

[2] Department of Mathematical Sciences, University of Bath, Bath, UK

proposed that the linear dispersion in the KdV equation with Fourier symbol $1 - \frac{1}{6}\xi^2$ should be replaced by the exact linear dispersion in the Euler equation with Fourier symbol

$$m(\xi) := \sqrt{\frac{\tanh(\xi)}{\xi}}.$$

Note that the dispersion in the KdV equation is the second-order approximation of $m$ at $\xi = 0$. This leads to the nonlinear nonlocal evolution equation

$$u_t + (K * u + u^2)_x = 0,$$

known as the *Whitham equation*. Here, $u(x, t)$ describes the one-dimensional wave profile and the integral kernel $K$ is given by

$$K(x) = (\mathcal{F}^{-1}m)(x) = \frac{1}{2\pi}\int_{\mathbb{R}} m(\xi) \exp(ix\xi) \, d\xi.$$

The function $K$ will be referred to as the *Whitham kernel*, and the function $m$ as the *Whitham symbol*. Specializing to traveling waves $u = \varphi(x - ct)$ where $c > 0$ is the wave speed, integrating and performing a Galilean change of variables, the Whitham equation reduces to the nonlinear integral equation

$$c\varphi - K * \varphi - \varphi^2 = 0. \tag{1}$$

We are interested in functions $\varphi\colon \mathbb{R} \to \mathbb{R}$ which satisfy (1) pointwise on $\mathbb{R}$, and which we refer to as solutions of (1) with wave speed $c$. More specifically, the results of this paper will concern *solitary* solutions, also called *solitary-wave* solutions. These are solutions $\varphi\colon \mathbb{R} \to \mathbb{R}$ satisfying $\lim_{|x|\to\infty}\varphi(x) = 0$.

Despite its simple form, the nonlocal and nonlinear nature of the Whitham equation has made it challenging to study. Recent years have seen a large amount of existence and qualitative results on the solutions of the equation. Traveling small-amplitude periodic solutions were found by Ehrnström and Kalisch [21] using the Crandall–Rabinowitz bifurcation theorem. Then, Ehrnström et al. [19] proved the existence of solitary waves using a variational method for a class of Whitham-type equations. This was followed up by Arnesen [4] where a class covering the Whitham equation with surface tension was considered. By applying a different technique—the implicit function theorem—Stefanov and Wright [31] achieved the same result. Ehrnström and Wahlén [22] showed the existence of a traveling cusped periodic wave $\varphi$ using global bifurcation theory, and proved that $\frac{c}{2} - \varphi(x) \sim |x|^{1/2}$ near the origin. This wave attains the highest amplitude possible and is referred to as an *extreme wave solution*. They also conjectured that $\varphi$ is convex and $\varphi = \frac{c}{2} - \sqrt{\pi/8}|x|^{1/2} + o(x)$ as $x \to 0$. Convexity of the extreme wave was shown by Encisco et al. [23] using a computer assisted proof.

The goal of this paper is to prove the existence of an extreme solitary-wave solution of (1) and our plan is to use a global bifurcation theorem appearing in [11]; see

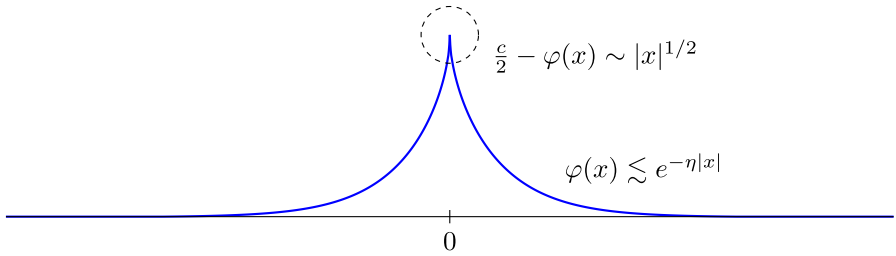

**Fig. 1** An extreme solitary-wave solution found by taking a limit of elements along the global bifurcation curve in Theorem 4.8. The wave speed $c$ is supercritical, that is, $c > 1$. The wave profile $\varphi$ is even and smooth on $\mathbb{R} \setminus \{0\}$. It has exponential decay as $|x| \to \infty$ and behaves like $\frac{c}{2} - C|x|^{1/2}$ as $|x| \to 0$. Ehrnström and Wahlén conjecture in [22] that $C = \sqrt{\pi/8}$

also [10,15,16]. The first main step is the construction of a local bifurcation curve, emanating from the point $(\varphi, c) = (0, 1)$, and the second is the application of the global bifurcation theorem.

A key to our success is the fact that a lot of qualitative properties have been shown for the Whitham kernel, the Whitham symbol and the solutions of (1), thanks to [9,21,22]. These guide us in choosing a convenient function space to study (1) and have been extremely useful in the application of the global bifurcation theorem. In Sect. 2, we list the relevant properties and prove an integral identity. We also study how sequences of solutions converge and the Fredholm properties of important linear operators.

Another key is the recently developed center manifold theorem for nonlocal equations in [24]. This result states that nonlocal equations with exponentially decaying convolution kernels are essentially local equations near an equilibrium. It also provides a method to derive the local equation, which can then be studied using familiar ODE tools. In our case, the equilibrium is $(\varphi, c) = (0, 1)$. Although the Whitham kernel has the required exponential decay, it fails a local integrability condition. Seeing that this condition is only for proving Fredholm properties of linear operators, we directly prove these properties instead. All necessary changes for the general center manifold theorem are listed in Appendix B. In Sect. 3, we state the center manifold theorem for the Whitham equation and compute the corresponding local equation. More specifically, we prove the following.

**Lemma 1.1** *There exist a neighborhood $\mathcal{V}' \subset \mathbb{R}$ of $c = 1$ and a number $\delta' > 0$ such that if $\varphi \colon \mathbb{R} \to \mathbb{R}$ satisfies $\sup_{y \in \mathbb{R}} \|\varphi(\cdot + y)\|_{H^3([0,1])} \leq \delta'$, then $\varphi$ solves (1) with wave speed $c \in \mathcal{V}'$ if and only if $\varphi$ solves the second-order ODE*

$$\varphi'' = -6\varphi^2 + \frac{19}{5}(\varphi')^2 + 6(c-1)\varphi + \mathcal{O}(|(\varphi, \varphi')|((c-1)^2 + |\varphi|^2 + |\varphi'|^2)). \quad (2)$$

The ODE in this lemma is a $c$-dependent family of perturbed KdV equations. Restricting to $c > 1$ with $c$ sufficiently close to 1, it features a unique positive even solitary-wave solution $\varphi$ with exponential decay for each fixed $c$. Using Eq. (2), we show that $\sup_{y \in \mathbb{R}} \|\varphi(\cdot + y)\|_{H^3([0,1])} \lesssim c - 1$. So for $c$ sufficiently close to 1, $\varphi$ is also a solution to equation (1). We thus arrive at the first main result of this paper (repeated as Theorem 3.3 in Sect. 3).

**Theorem 1.2** *There exists a unique local bifurcation curve $\mathcal{C}_{loc}$ which emanates from $(\varphi, c) = (0, 1)$ and consists of the non-trivial even solitary-wave solutions $\varphi$ to* (1) *with wave speeds $c \in \mathcal{V}'$ satisfying $\sup_{y \in \mathbb{R}} \|\varphi(\cdot + y)\|_{H^3([0,1])} < \delta'$.*

While both [19] and [31] contain existence results for supercritical solitary waves, the additional information provided by the center manifold approach concerning uniqueness is crucial in the subsequent analysis. To end Sect. 3, we use the center manifold theorem to prove that the linearization of the left-hand side of (1) along $\mathcal{C}_{loc}$ is invertible. This is in preparation for the global bifurcation theorem.

The global bifurcation theory in [11] can now be applied to extend $\mathcal{C}_{loc}$ and this extension is referred to as the global bifurcation curve $\mathcal{C}$. The theory dictates several possible behaviors for $\mathcal{C}$ and the content of Sect. 4 is the exclusion of unwanted behaviors. We rule out the loss of compactness alternative using qualitative properties of the solutions, how sequences of solutions converge and an integral identity for (1). Then, we show that the blowup alternative happens as the Sobolev norm blows up and that an extreme solitary-wave solution is obtained in the limit. More precisely, we have the following result (repeated as Theorem 4.8 in Sect. 4).

**Theorem 1.3** *There exists a sequence of elements $(\varphi_n, c_n)$ on the global bifurcation curve $\mathcal{C}$ such that $\lim_{n \to \infty} \|\varphi_n\|_{H^3} = \infty$ and $(\varphi_n, c_n) \to (\varphi, c)$ locally uniformly, where $\varphi$ is a solitary-wave solution of (1) with supercritical wave speed $c > 1$. The solitary solution $\varphi$ is even, bounded, continuous, exponentially decaying, smooth everywhere except at the origin and*

$$C_1 |x|^{1/2} \leq \frac{c}{2} - \varphi(x) \leq C_2 |x|^{1/2} \text{ as } |x| \to 0,$$

*for some constants $0 < C_1 < C_2$.*

The function $\varphi$ in the above theorem is the extreme solitary-wave solution we set out to find and is illustrated in Fig. 1.

By demonstrating the use of recent spatial-dynamics tools, this paper serves as an example to studies of other nonlocal nonlinear evolution equations. In particular, these results will likely extend to a larger class of equation, such as in [5] and [20].

Finally, it is interesting to compare our results with the global bifurcation theory for the water wave problem. The existence of an unbounded, connected set of solitary water waves, including a highest wave in a certain limit, was proved by Amick and Toland [2,3] following several earlier small-amplitude results. Around the same time, Amick et al. [1] verified Stokes' conjecture for both periodic and solitary water waves, showing in particular that the limiting solitary wave is Lipschitz continuous at the crest with a corner enclosing a 120° angle. Thus, the behavior at the crest is different from the extreme Whitham wave, which has no corner due to the $C^{1/2}$ cusp in Theorem 1.3. The construction of the global solution continua in [2,3] is also different from ours. While both proofs are based on nonlinear integral equations, the common approach in [2,3] is to first apply global bifurcation theory to a regularized problem and then pass to the limit. On the other hand, we use global bifurcation theory directly on the solitary Whitham problem. A similar approach has in fact recently been used for solitary water waves with vorticity and stratification, but based on a PDE formulation [11,32]. For the water wave problem with vorticity and stratification, the limiting behavior of large-amplitude waves is more complex and there is numerical and some analytical evidence of overhanging waves; see for example [14,17,18,28] and references therein.

**Notation**

We use the following notations for function spaces.

- The space of $p$th power integrable functions on an interval $I \subset \mathbb{R}$ with respect to a measure $\mu$ is denoted by

$$L^p(I, \mu) := \left\{ f : \mathbb{R} \to \mathbb{R} \,\Big|\, \|f\|_{L^p(\mu)} < \infty \right\},$$

where

$$\|f\|_{L^p(I,\mu)} := \left( \int_I |f|^p \, d\mu \right)^{1/p} \quad \text{if } p \in [1, \infty)$$

and

$$\|f\|_{L^\infty(I,\mu)} := \mu\text{-ess-sup}_{x \in I} |f(x)| \quad \text{if } p = \infty.$$

For $\sigma \in \mathbb{R}$, we write

$$L^p(I)| := L^p(I, dx), \quad L^p := L^p(\mathbb{R}, dx) \quad \text{and} \quad L^p_\sigma := L^p(\mathbb{R}, \omega^p_\sigma \cdot dx),$$

where $dx$ is the Lebesgue measure and $\omega_\sigma : \mathbb{R} \to \mathbb{R}$ is a positive and smooth function, which equals $\exp(\sigma|x|)$ for $|x| > 1$. In particular, functions in $L^p_\eta$ when $\eta > 0$ are necessarily exponentially decaying while functions in $L^p_{-\eta}$ can grow exponentially.

- The Sobolev space is denoted by

$$W^{j,p}(I) := \left\{ f : I \to \mathbb{R} \,\Big|\, f^{(n)} \in L^p(I), \text{ for } 0 \leq n \leq j \right\}$$

and the weighted Sobolev space is

$$W^{j,p}_\sigma := \left\{ f : \mathbb{R} \to \mathbb{R} \,\Big|\, f^{(n)} \in L^p_\sigma, \text{ for } 0 \leq n \leq j \right\},$$

where $f^{(n)}$ are weak derivatives of $f$ for $1 \leq n \leq j$. These spaces are equipped with the norms

$$\|f\|_{W^{j,p}} := \left( \sum_{n=0}^j \|f^{(n)}\|^p_{L^p} \right)^{1/p} \quad \text{and} \quad \|f\|_{W^{j,p}_\sigma} := \left( \sum_{n=0}^j \|f^{(n)}\|^p_{L^p_\sigma} \right)^{1/p}.$$

We have the natural inclusions $W^{j,p}_{\sigma_2} \subset W^{j,p}_{\sigma_1}$ whenever $\sigma_1 < \sigma_2$. For $p = 2$, we denote the Hilbert spaces $W^{j,2}(I)$ and $W^{j,2}_{\sigma}$ by $H^j(I)$ and $H^j_{\sigma}$, respectively. As before, when $I = \mathbb{R}$, we omit writing $\mathbb{R}$.

- We define the space of uniformly local $H^j$ functions

$$H^j_u := \left\{ f \in H^j_{\text{loc}} \,\middle|\, \|f\|_{H^j_u} < \infty \right\}, \quad \text{where} \quad \|f\|_{H^j_u} := \sup_{y \in \mathbb{R}} \|f(\cdot + y)\|_{H^j([0,1])}.$$

- $C^k$ denotes the space of $k$ times continuously differentiable functions $f \colon \mathbb{R} \to \mathbb{R}$. $BUC^k \subset C^k$ denotes the space of functions with bounded and uniformly continuous derivatives of order up to and including $k$. $C^{k,\alpha}$ denotes the Hölder spaces

$$C^{k,\alpha} := \left\{ f \in BUC^k \,\middle|\, \sup_{h \neq 0} \frac{|f^{(k)}(x+h) - f^{(k)}(x)|}{|h|^{\alpha}} < \infty \right\}.$$

- The Besov space is denoted by $B^s_{p,q}$, where $s \in \mathbb{R}$, $p, q \in [1, \infty]$. We have

$$B^s_{2,2} = H^s, \ s \in \mathbb{R} \quad \text{and} \quad B^s_{\infty,\infty} = C^{\lfloor s \rfloor, s - \lfloor s \rfloor}, \ s \in \mathbb{R}_+ \setminus \mathbb{N}.$$

- $\mathscr{C}^k(\mathcal{X}, \mathcal{Y})$ denotes the space of $k$ times Fréchet differentiable mappings between two normed spaces $\mathcal{X}$ and $\mathcal{Y}$.

We use the following scaling of the Fourier transform:

$$\mathcal{F}f(\xi) := \int_{\mathbb{R}} f(x) \exp(-ix\xi) \, dx$$

and

$$\mathcal{F}^{-1}g(x) := \frac{1}{2\pi} \mathcal{F}g(-x).$$

## 2 Qualitative properties

### 2.1 The Whitham kernel and the Whitham symbol

The Whitham kernel $K$ is given by $(\mathcal{F}^{-1}m)(x)$, where

$$m(\xi) = \sqrt{\frac{\tanh(\xi)}{\xi}} = 1 - \frac{1}{6}\xi^2 + \frac{19}{360}\xi^4 + \mathcal{O}(\xi^6). \tag{3}$$

Since $m(0) = 1$, we have $\int_{\mathbb{R}} K \, dx = 1$. However, since $m \notin L^1$, $K$ is singular at the origin. More specifically,

$$K(x) = \frac{1}{\sqrt{2\pi|x|}} + K_{reg}(x), \qquad (4)$$

where $K_{reg}$ is real analytic on $\mathbb{R}$; see Proposition 2.4 in [22]. In addition, as $|x| \to \infty$,

$$K(x) = \frac{\sqrt{2}}{\pi\sqrt{|x|}} \exp\left(-\frac{\pi}{2}|x|\right) + \mathcal{O}\left(|x|^{-3/2} \exp\left(-\frac{\pi}{2}|x|\right)\right), \qquad (5)$$

by Corollary 2.26 in [22]. Since $m$ is an even function, so is $K$. The fact that $K$ is a positive function has been shown in Proposition 2.23 in [22].

Because the Whitham symbol $m$ satisfies

$$|m^{(j)}(\xi)| \leq C_j (1 + |\xi|)^{-\frac{1}{2} - j}, \quad \text{for } j \in \mathbb{N},$$

the linear operator

$$\varphi \mapsto m(D)\varphi := K * \varphi, \quad B_{p,q}^s \to B_{p,q}^{s+1/2} \qquad (6)$$

is bounded; see for example Proposition 2.78 in [6].

## 2.2 Properties of solutions

When choosing appropriate function spaces for (1), we will rely on the following qualitative properties of solutions.

**Proposition 2.1** *Let $\varphi$ be a continuous and bounded solution to (1) with wave speed $c \geq 1$. We have*

(i) *(non-negativity) $\varphi \geq 0$;*
(ii) *(exponential decay) if $\varphi$ is solitary and $c > 1$, there exists $\eta > 0$ such that $\exp(\eta|x|)\varphi(x) \in L^\infty(\mathbb{R})$, that is, $\varphi$ has exponential decay;*
(iii) *(symmetry) if $\varphi$ is solitary, $c > 1$ and $\sup_{x \in \mathbb{R}} \varphi(x) < c/2$, there exists $x_0 \in \mathbb{R}$ such that $\varphi(\cdot - x_0)$ is an even function which is non-increasing on $[0, \infty)$;*
(iv) *(regularity) $\varphi$ is smooth on any open set where $\varphi < c/2$;*
(v) *(boundedness) if $\varphi$ is of class $BUC^1$, even, non-constant and non-increasing on $(0, \infty)$, then $\varphi' < 0$ and $\varphi < c/2$ on $(0, \infty)$. If $\varphi$ in addition is of class $BUC^2$, then $\varphi < c/2$ everywhere;*
(vi) *(singularity) if $\varphi$ is even, non-constant, non-increasing on $(0, \infty)$, $\sup_{x \in \mathbb{R}} \varphi(x) \leq c/2$ and $\varphi(0) = c/2$, then as $|x| \to 0$,*

$$C_1 |x|^{\frac{c}{2}} \leq \frac{1+c}{2} - \varphi(x) \leq C_2 |x|^{\frac{1}{2}};$$

(vii) *(lower bound on the wave speed)* if $\varphi$ is non-trivial and has finite limits $\lim_{x \to \pm\infty} \varphi(x)$, then $c > 1$;

(viii) *(upper bound on the wave speed)* if $\varphi$ is non-constant and $\varphi \leq c/2$, then $c \leq 2$.

Item (i) is stated in Lemma 4.1 in [22]. Items (ii) and (iii) are Proposition 3.13 and Theorem 4.4 in [9] respectively. The optimal exponent $\eta = \eta_c$ depends on $c$ and is given implicitly by $\sqrt{\tan(\eta_c)/\eta_c} = c$, with $\eta_c \in (0, \pi/2)$; see [5], Theorem 6.2. We remark that the requirement $\sup_{x \in \mathbb{R}} \varphi(x) < c/2$ in (iii) is not mentioned in [9] despite its importance in the proof; see the introduction in [30] for a detailed discussion. Items (iv), (v) and (vi) can be found in Theorem 4.9, Theorem 5.1 and Theorem 5.4 in [22]. The upper bound in (viii) comes from [22], Eq. (6.9). The lower bound in (vii) comes from non-negativity and the following proposition.

**Proposition 2.2** *Let $\varphi$ be a bounded and continuously differentiable solution to* (1) *with wave speed $c$, such that the limits $\lim_{x \to \pm\infty} \varphi(x)$ exist. Then*

$$\lim_{R \to \infty} \int_{-R}^{R} \varphi(\varphi - (c - 1)) \, dx = 0. \tag{7}$$

*In particular,* $\inf \varphi < c - 1 < \sup \varphi$, $\varphi \equiv c - 1$ *or* $\varphi \equiv 0$ *if* $c \geq 1$. *The latter statement is in fact true for bounded and continuous solutions $\varphi$ with finite limits $\lim_{x \to \pm\infty} \varphi(x)$.*

**Proof** In general, if $\varphi$ is any bounded and continuously differentiable function, the limits $\lim_{x \to \pm\infty} \varphi(x)$ exist, and $\mathcal{K}$ is any non-negative even function with $\int_{\mathbb{R}} \mathcal{K} \, dx = 1$, then

$$\lim_{R \to \infty} \int_{-R}^{R} (\varphi - \mathcal{K} * \varphi) \, dx = 0.$$

A proof of this can be found in [8], pp. 113–114.

Since $\varphi$ solves (1) with wave speed $c$,

$$K * \varphi = c\varphi - \varphi^2,$$

and since the Whitham kernel $K$ is a non-negative even function with $\int_{\mathbb{R}} K \, dx = 1$,

$$\begin{aligned} 0 &= \lim_{R \to \infty} \int_{-R}^{R} (\varphi - K * \varphi) \, dx \\ &= \lim_{R \to \infty} \int_{-R}^{R} (\varphi - c\varphi + \varphi^2) \, dx \\ &= \lim_{R \to \infty} \int_{-R}^{R} \varphi(\varphi - (c - 1)) \, dx. \end{aligned}$$

By non-negativity of bounded and continuous solutions with $c \geq 1$, $\varphi - (c-1)$ must be sign-changing, $\varphi \equiv 0$, or $\varphi \equiv c - 1$, otherwise the generalized integral cannot converge to zero. This proves the claim for bounded and continuously differentiable functions $\varphi$.

If $\varphi$ is only a bounded and continuous solution, convolution with a non-zero smooth and compactly supported test function $\phi \geq 0$ gives

$$\phi * (\varphi - K * \varphi) = \phi * \varphi(\varphi - (c-1)).$$

The left-hand side equals

$$\phi * \varphi - K * (\phi * \varphi),$$

due to associativity and commutativity of convolution. By Lebesgue's dominated convergence theorem, the function $\phi * \varphi$ is bounded, continuously differentiable and the limits as $x \to \pm\infty$ exist. It follows that

$$\lim_{R \to \infty} \int_{-R}^{R} [\phi * \varphi - K * (\phi * \varphi)] \, dx = 0,$$

which implies

$$\lim_{R \to \infty} \int_{-R}^{R} \phi * (\varphi(\varphi - (c-1))) \, dx = 0.$$

Again, we must have $\varphi \equiv 0$, $\varphi \equiv c - 1$, or $\varphi - (c-1)$ is sign-changing. □

### 2.3 Convergence of solution sequences

Modes of convergence of solution sequences will be important in ruling out alternatives from the global bifurcation theorem. We start with pointwise convergence, using the Arzelà—Ascoli theorem and the smoothing property of convolution with $K$.

**Proposition 2.3** *Let $(\varphi_n)_{n=1}^{\infty}$ be a sequence of continuous and bounded solutions to* (1) *such that each $\varphi_n$ has wave speed $c_n \in [1, 2]$ and $\sup_{x \in \mathbb{R}} \varphi_n(x) \leq c_n/2$. Then, there exists a subsequence $(\varphi_{n_k})_{k=1}^{\infty}$ satisfying*

$$\lim_{k \to \infty} c_{n_k} = c \in [1, 2], \quad \lim_{k \to \infty} \varphi_{n_k}(x) = \varphi(x),$$

*for every $x \in \mathbb{R}$. The convergence is uniform on every bounded interval of $\mathbb{R}$. The limit $\varphi$ is a continuous, bounded and non-negative solution of* (1) *with wave speed $c$, and $\sup_{x \in \mathbb{R}} \varphi(x) \leq c/2$.*

**Proof** We can without loss of generality assume that $\lim_{n\to\infty} c_n = c \in [1, 2]$. For each $n$, we have

$$\left(\frac{c_n}{2} - \varphi_n\right)^2 - \frac{c_n^2}{4} = -K * \varphi_n.$$

A rearrangement gives

$$\varphi_n = \frac{c_n}{2} - \sqrt{\left(\frac{c_n}{2}\right)^2 - K * \varphi_n}.$$

We claim that the right-hand side forms an equicontinuous sequence. Indeed, $\varphi \mapsto K * \varphi$ is a bounded map from $L^\infty \subset B^0_{\infty,\infty}$ to $B^{1/2}_{\infty,\infty} = C^{1/2}$ according to (6). Because $c_n \in [1, 2]$, this gives

$$\|K * \varphi_n\|_{C^{1/2}} \lesssim \sup_{x \in \mathbb{R}} \varphi_n(x) \leq \frac{c_n}{2} \leq 1.$$

Hence, $(K * \varphi_n)_{n=1}^\infty$ is an equicontinuous sequence of functions. The square root of a non-negative equicontinuous sequence is an equicontinuous sequence. So, $(\varphi_n)_{n=1}^\infty$ is equicontinuous. The Arzelà—Ascoli theorem gives a subsequence $(\varphi_{n_k})_{k=1}^\infty$, which converges uniformly to a function $\varphi$ on each bounded interval of $\mathbb{R}$. Also, $\varphi$ is continuous and bounded by $c/2$.

Finally, since $\sup_{x \in \mathbb{R}} \varphi_n(x) \leq 1$ and $\|K\|_{L^1} = 1$, Lebesgue's dominated convergence theorem gives $K * \varphi_n(x) \to K * \varphi(x)$ as $n \to \infty$ for all $x$. It follows that $\varphi$ is a solution to (1) with wave speed $c$. □

Here are several immediate consequences.

**Corollary 2.4**

(i) *If $\varphi$ solves (1) with wave speed $c \in [1, 2]$, $\sup_{x \in \mathbb{R}} \varphi(x) \leq c/2$, and $\lim_{x \to \infty} \varphi(x) = a$, then $a$ solves (1) with wave speed $c$. In particular, the constant $a$ is either zero or $c - 1$.*

(ii) *Let $\varphi_n$, $c_n$, $\varphi$ and $c$ be as in Proposition 2.3. If $\varphi_n$ is even and monotone on $[0, \infty)$ for each $n$, its locally uniform limit $\varphi$ inherits evenness and monotonicity on $[0, \infty)$. If in addition $\lim_{|x| \to \infty} \varphi(x) = 0$, then $\varphi_n$ converges to $\varphi$ uniformly on $\mathbb{R}$.*

(iii) *Let $(\varphi_n)_{n=1}^\infty$ be a sequence of even solutions which are decreasing on $[0, \infty)$. Define $\tau_{x_n} \varphi_n := \varphi_n(\cdot + x_n)$ for a sequence of real numbers $x_n$ with $\lim_{n \to \infty} x_n = \infty$. Then, the sequence of translated solutions $\tau_{x_n} \varphi_n$ is a sequence of solutions to (1). It has a non-increasing locally uniform limit $\tilde{\varphi}$.*

**Proof** Item (iii) is straightforward. We only prove items (i) and (ii).

For (i), we define

$$\tau_n \varphi := \varphi(\cdot + n), \quad n \in \mathbb{N}.$$

Each $\tau_n \varphi$ is a solution of (1) with wave speed $c \in [1, 2]$. We have

$$\lim_{n \to \infty} \tau_n \varphi(x) = \lim_{n \to \infty} \varphi(x + n) = \lim_{x \to \infty} \varphi(x) = a, \quad \text{for each } x \in \mathbb{R}.$$

By Proposition 2.3, $a$ is a constant solution to (1) with wave speed $c$ and hence by Proposition 2.2, we have $a = 0$ or $a = c - 1$.

The evenness and monotonicity in (ii) are clear, so assume that $\lim_{|x| \to \infty} \varphi(x) = 0$ and fix $\epsilon > 0$. Then, there exists $R > 0$ such that

$$|\varphi(x)| \leq \epsilon \text{ for } |x| \geq R. \tag{8}$$

Due to $\lim_{n \to \infty} \varphi_n(x) = \varphi(x)$ locally uniformly, there exists an $N_\epsilon > 0$ such that

$$\sup_{|x| \leq R} |\varphi_n(x) - \varphi(x)| \leq \epsilon \text{ for } n > N_\epsilon.$$

In particular, we have $|\varphi_n(R) - \varphi(R)| \leq \epsilon$, which in turn implies $|\varphi_n(R)| \leq 2\epsilon$ by (8). Since $\varphi_n$ is non-increasing, we have $|\varphi_n(x)| \leq 2\epsilon$ for $|x| \geq R$. But then, again by (8), $|\varphi_n(x) - \varphi(x)| \leq 3\epsilon$ for $|x| \geq R$ and $n > N_\epsilon$, and the claim about uniform convergence is proved. □

Now, we consider convergence in $H^j$ for $j > 0$. Combining the smoothing property of convolution with $K$ and (1), we use a bootstrap argument to increase the regularity of the solutions, starting with convergence in $L^2 = H^0$.

**Proposition 2.5** *Let $\varphi_n, c_n, \varphi$ and $c$ be as in Proposition 2.3. If*

(a) $\varphi_n \to \varphi$ *uniformly and* $\varphi_n \to \varphi$ *in* $L^2$,
(b) $\sup_{x \in \mathbb{R}} \varphi_n(x) < c_n/2$ *and* $\sup_{x \in \mathbb{R}} \varphi(x) < c/2$,

*then $\varphi_n \to \varphi$ in $H^j$ for any $j > 0$.*

**Proof** Since $\varphi_n$ and $\varphi$ solve (1) with wave speed $c_n$ and $c$ respectively, we can write

$$\varphi_n - \varphi = f(K * \varphi_n, c_n) - f(K * \varphi, c), \quad \text{where} \quad f(\omega, c) = \frac{c}{2} - \sqrt{\frac{c^2}{4} - \omega}. \tag{9}$$

Letting $\omega_n = K * \varphi_n$ and $\omega = K * \varphi$, the assumptions imply

(a') $\omega_n \to \omega$ in $H^{1/2}$ and $\omega_n \to \omega$ uniformly,
(b') $\inf_{n,x} \left( \frac{c_n^2}{4} - \omega_n \right) > \varepsilon/3$,
(c') $\inf_{n,x} \left( \frac{c^2}{4} - \omega_n \right) > \varepsilon/3$,

for some $\varepsilon > 0$ and for sufficiently large $n$. Without loss of generality, (b') and (c') are assumed to hold for all $n$. We claim that if (a'), (b') and (c') are met, then

$$f(\omega_n, c_n) \to f(\omega, c) \text{ in } H^j, \quad \text{for some } j \in (0, 1/2).$$

Consider

$$\|f(\omega_n, c_n) - f(\omega, c)\|_{H^j} \leq \|f(\omega_n, c_n) - f(\omega_n, c)\|_{H^j} + \|f(\omega_n, c) - f(\omega, c)\|_{H^j}.$$

A quick calculation gives

$$f(\omega_n, c_n) - f(\omega_n, c) = \frac{c_n - c}{2} \cdot \left(1 - \frac{1}{2} \cdot \frac{c_n + c}{\sqrt{\frac{c_n^2}{4} - \omega_n} + \sqrt{\frac{c^2}{4} - \omega_n}}\right).$$

Define

$$g_n(x) = 1 - \frac{1}{2} \cdot \frac{c_n + c}{\sqrt{\frac{c_n^2}{4} - x} + \sqrt{\frac{c^2}{4} - x}}, \quad \text{with domain } D_{g_n} = \left[0, \min\left\{\frac{c^2}{4}, \frac{c_n^2}{4}\right\} - \frac{\varepsilon}{3}\right).$$

For each $n$, $g_n$ is smooth and $g_n(0) = 0$. Moreover, the range of $\omega_n$ belongs to the domain of $g_n$. A standard result in the theory of paradifferential operators, for instance Theorem 2.87 in [6], gives

$$\|g_n(\omega_n)\|_{H^j} \leq C \|\omega_n\|_{H^j},$$

where $C > 0$ depends on $j$, $\sup_{x \in \mathbb{R}} |\omega_n|$ and $\sup_{x \in D_{g_n}} |g_n'(x)|$. A computation shows that $|g_n'|$ is uniformly bounded in $n$ for $c_n, c \in [1, 2]$ and $x \in D_{g_n}$. Also, $\sup_n \|\omega_n\|_{H^{1/2}} < \infty$, as well as $\sup_{n,x} |\omega_n(x)| < \infty$ by (a'). It follows that $\|g_n(\omega_n)\|_{H^j}$ is uniformly bounded in $n$ and

$$\|f(\omega_n, c_n) - f(\omega_n, c)\|_{H^j} = \frac{1}{2}|c_n - c| \cdot \|g_n(\omega_n)\|_{H^j} \to 0, \quad \text{as } n \to \infty.$$

To deal with the second term, let $c$ be fixed and define

$$h(x) = f'(x, c) - f'(0, c), \quad \text{with domain } D_h = \left[0, \frac{c^2}{4} - \frac{\varepsilon}{3}\right).$$

Then, we can write

$$f(\omega_n, c) - f(\omega, c) = (\omega_n - \omega) \int_0^1 h(\omega_n + \tau(\omega - \omega_n)) \, d\tau + f'(0, c)(\omega_n - \omega).$$

Due to (a') and (b'), we have $\omega_n(s) + \tau(\omega(s) - \omega_n(s)) \in D_h$ for all $s \in \mathbb{R}$. Note that $h$ is smooth and $h(0) = 0$. Theorem 2.87 in [6] gives an estimate for the integral

$$\left\| \int_0^1 h(\omega_n + \tau(\omega - \omega_n)) \, d\tau \right\|_{H^{1/2}} \leq \sup_{\tau \in [0,1]} \| h(\omega_n + \tau(\omega - \omega_n)) \|_{H^{1/2}}$$
$$\leq C \sup_{\tau \in [0,1]} \| \omega_n + \tau(\omega - \omega_n) \|_{H^{1/2}},$$

where $C$ depends on $s$, $\sup_x |f''(x, c)|$, $c$ and $\varepsilon$. Another standard result in paradifferential calculus, for example Theorem 8.3.1 in [29], gives

$$\| f(\omega_n, c) - f(\omega, c) \|_{H^j} \lesssim \| \omega_n - \omega \|_{H^{1/2}} \left\| \int_0^1 h(\omega_n + \tau(\omega - \omega_n)) \, d\tau \right\|_{H^{1/2}}$$
$$+ |f'(0, c)| \cdot \| \omega_n - \omega \|_{H^j},$$

for all $j \in (0, 1/2)$, which tends to 0 as $n \to \infty$ by (a').

So, the right-hand side of (9) tends to 0 in $H^j$ for all $j \in (0, 1/2)$. Hence, $\varphi_n \to \varphi$ in $H^j$. Then, convolution with $K$ increases the regularity of $\varphi_n$ and $\varphi$ by $1/2$. Choosing $j = 1/4$ and replacing (a') with

$$\omega_n \to \omega \quad \text{in } H^{1/4+1/2},$$

the critical case of Theorem 8.3.1 in [29] is no longer relevant and the convergence is in $H^{j'}$ for all $j' \in (0, 1/4 + 1/2]$. By iterating as many times as needed, the claim of the proposition is proved. □

**Remark 2.6** Since Theorem 2.87 in [6] and Theorem 8.3.1 in [29] are valid for the Besov spaces $B_{p,q}^s$, we can replace the $L^2$ space in (a) with $B_{p,q}^0$ and obtain $\varphi_n \to \varphi$ in $B_{p,q}^s$ for $s > 0$ using the same proof idea.

### 2.4 Fredholmness of linear operators

Let $j > 0$ be an integer. As a preparation for future bifurcation results, we study the operators

$$\mathcal{T} \colon \varphi \mapsto \varphi - K * \varphi, \quad H_{-\eta}^j \to H_{-\eta}^j$$

and

$$L[\varphi^*, c^*] \colon \phi \mapsto c^* \phi - K * \phi - 2\varphi^* \phi, \quad H^j \to H^j,$$

where $\varphi^*$ is a solitary-wave solution with wave speed $c^* > 1$, satisfying $\sup_{x \in \mathbb{R}} \varphi^*(x) < c^*/2$. These are linearizations of the left-hand side in (1) at $(0, 1)$ and $(\varphi^*, c^*)$, respectively. We show that $\mathcal{T}$ and $L[\varphi^*, c^*]$ are Fredholm with Fredholm index two and zero respectively, using results from [27]. The central idea is to relate a pseudodifferential operator $t(x, D) \colon H^j \to H^j$ to a positively homogeneous function $A$ via the symbol $t(x, \xi)$. By studying the boundary value and the winding number of $A$ around the origin, the Fredholm property of $t(x, D)$ can be determined. Appendix A summarizes the relevant theorems from [27].

Up until now, the weight $\eta > 0$ has remained somewhat mysterious. Since our interest lies in the Fredholm properties of $\mathcal{T}$, $\eta$ should be chosen so that $\mathcal{T}$ is at least bounded. By (4), $K$ is locally $L^1$ around the origin. From (5), we deduce that $K \in L^1_{\eta_0}$ for $\eta_0 \in (0, \pi/2)$. It follows that $\mathcal{T} \colon H^j_{-\eta} \to H^j_{-\eta}$ is bounded for any $\eta \in (0, \pi/2)$. Indeed, since $\omega_{-\eta}$ is 1 on $[-1, 1]$ and equals $\exp(-\eta|x|)$ as $|x| \to \infty$, we have

$$\|K * \varphi\|^2_{L^2_{-\eta}} \lesssim \int_{\mathbb{R}} \left( \int_{\mathbb{R}} K(y) \varphi(x-y) \, \mathrm{d}y \right)^2 \exp(-2\eta|x|) \, \mathrm{d}x$$

$$\leq \int_{\mathbb{R}} \left( \int_{\mathbb{R}} K(y) \varphi(x-y) \, \mathrm{d}y \right)^2 \exp\left( -2\eta|x-y| + 2\eta|y| \right) \, \mathrm{d}x$$

$$= \int_{\mathbb{R}} \left( \int_{\mathbb{R}} K(y) \exp(\eta|y|) \cdot \varphi(x-y) \exp(-\eta|x-y|) \, \mathrm{d}y \right)^2 \, \mathrm{d}x$$

$$= \left\| \left[ K \cdot \exp(\eta|\cdot|) \right] * \left[ \varphi \cdot \exp(-\eta|\cdot|) \right] \right\|^2_{L^2}$$

$$\leq \|K\|^2_{L^1_\eta} \|\varphi\|^2_{L^2_{-\eta}},$$

where Young's inequality for the $L^p$ norms of convolutions is used in the last step. Noting that

$$\frac{\mathrm{d}^n}{\mathrm{d}x^n}(K * \varphi) = K * \varphi^{(n)}$$

and applying the above estimates on $\varphi^{(n)}$,

$$\|\mathcal{T}\varphi\|_{H^j_{-\eta}} \lesssim \|\varphi\|_{H^j_{-\eta}},$$

for $\eta$ in the range $(0, \pi/2)$.

### 2.4.1 Fredholmness of $\mathcal{T}$

Multiplication with $\cosh(\eta \cdot)$

$$M_{\cosh} \colon \varphi \mapsto \cosh(\eta \cdot) \varphi, \quad H^j \to H^j_{-\eta}.$$

is an invertible linear operator. Its inverse is multiplication with $1/\cosh(\eta \cdot)$, mapping $H^j_{-\eta}$ to $H^j$. Conjugating $\mathcal{T}$ with these gives

$$\tilde{\mathcal{T}} = M_{\cosh}^{-1} \circ \mathcal{T} \circ M_{\cosh}, \quad H^j \to H^j,$$

and more explicitly

$$\tilde{\mathcal{T}}\varphi(x) = \mathrm{Id} - (\tilde{K}(\,\cdot\,, x) * \varphi)(x), \text{ where } \tilde{K}(z, x) = K(z)\frac{\cosh(\eta(x-z))}{\cosh(\eta x)}.$$

Setting

$$\phi_\pm(x) = \frac{\exp(\pm \eta x)}{2\cosh(\eta x)} \text{ and } K_\pm(z) = K(z)\exp(\pm \eta z),$$

$\tilde{\mathcal{T}}$ can be rewritten as

$$\tilde{\mathcal{T}}\varphi(x) = \varphi(x) - \phi_+(x)(K_- * \varphi)(x) - \phi_-(x)(K_+ * \varphi)(x),$$

with the symbol

$$\tilde{t}(x, \xi) = 1 - \phi_+(x)m(\xi - i\eta) - \phi_-(x)m(\xi + i\eta).$$

**Lemma 2.7** *The conjugated pseudodifferential operator $\tilde{\mathcal{T}} = \tilde{t}(x, \mathrm{D}) \colon H^j \to H^j$ is a Fredholm operator.*

**Proof** The idea is to apply Proposition A.1. We define a positively homogeneous function $A$ by

$$A(x_0, x, \xi_0, \xi) := \tilde{t}\left(\frac{x}{x_0}, \frac{\xi}{\xi_0}\right)$$

for $x, \xi \in \mathbb{R}$ and $x_0, \xi_0 > 0$. In order to apply Proposition A.1, we need to check that $A$ is smooth in

$$\overline{\mathbb{S}} := \{(x_0, x, \xi_0, \xi) \in \mathbb{R}^4 \mid x_0^2 + x^2 = \xi_0^2 + \xi^2 = 1, x_0 \geq 0, \xi_0 \geq 0\},$$

and that $A(x_0, x, \xi_0, \xi) \neq 0$ on $\Gamma$, where $\Gamma$ is the boundary of $\overline{\mathbb{S}}$. $\Gamma$ can be decomposed into the arcs

$$\begin{aligned}
\Gamma_1 &= \{(0, 1, \xi_0, \xi) \mid \xi_0^2 + \xi^2 = 1, \xi_0 \geq 0\}, \\
\Gamma_2 &= \{(0, -1, \xi_0, \xi) \mid \xi_0^2 + \xi^2 = 1, \xi_0 \geq 0\}, \\
\Gamma_3 &= \{(x_0, x, 0, 1) \mid x_0^2 + x^2 = 1, x_0 \geq 0\}, \\
\Gamma_4 &= \{(x_0, x, 0, -1) \mid x_0^2 + x^2 = 1, x_0 \geq 0\}.
\end{aligned}$$

We compute the value of $A$ along each arc $\Gamma_i$ and show that $A$ is nowhere vanishing. From the computations, it will be apparent that $A$ is smooth on $\overline{\mathbb{S}}$.

On $\Gamma_1^* := \Gamma_1 \setminus \{(0, 1, 0, 1), (0, 1, 0, -1)\}$, we have $\xi_0 = \sqrt{1 - \xi^2} > 0$ and

$$A(x_0, x, \xi_0, \xi)\Big|_{\Gamma_1^*} = \lim_{x_0 \to 0^+} \tilde{t}\left(\frac{1}{x_0}, \frac{\xi}{\sqrt{1-\xi^2}}\right)$$

$$= 1 - m\left(\frac{\xi}{\sqrt{1-\xi^2}} - i\eta\right),$$

since $\lim_{y \to \infty} \phi_+(y) = 1$ and $\lim_{y \to \infty} \phi_-(y) = 0$. Let $\theta = \xi(1-\xi^2)^{-1/2}$. As $\xi \in (-1, 1)$, $\theta \in (-\infty, \infty)$. To compute the values at the endpoint $(0, 1, 0, 1)$, we note that taking the limit as $\xi \to 1^-$ corresponds to taking the limit as $\theta \to \infty$. A calculation shows that

$$|m(\theta \pm i\eta)|^4 = \frac{\sinh^2(2\theta) + \sin^2(2\eta)}{(\theta^2 + \eta^2)(\cosh(2\theta) + \cos(2\eta))^2}$$

which gives

$$\lim_{\theta \to \pm\infty} |m(\theta \pm i\eta)| = 0.$$

So, the value of $A$ at $(0, 1, 0, 1)$ is 1. Similarly, $A$ at $(0, 1, 0, -1)$ corresponds to taking the limit as $\theta \to -\infty$ and the value is 1. Along $\Gamma_1^*$, if $1 - m(\theta - i\eta) = 0$ for some $\theta \in \mathbb{R}$, then

$$\operatorname{Re}(m(\theta - i\eta)^2) = 1 \text{ and } \operatorname{Im}(m(\theta - i\eta)^2) = 0,$$

where

$$\operatorname{Re}(m(\theta - i\eta)^2) = \frac{\theta \sinh(2\theta) + \eta \sin(2\eta)}{(\theta^2 + \eta^2)(\cosh(2\theta) + \cos(2\eta))}$$

and

$$\operatorname{Im}(m(\theta - i\eta)^2) = \frac{\eta \sinh(2\theta) - \theta \sin(2\eta)}{(\theta^2 + \eta^2)(\cosh(2\theta) + \cos(2\eta))}.$$

The second equation is satisfied only when the numerator $\eta \sinh(2\theta) - \theta \sin(2\eta)$ is zero. Since this is an odd and smooth function in $\theta$, a trivial solution is $\theta = 0$. Since

$$\frac{\sinh(2\theta)}{2\theta} > 1 \text{ and } \frac{\sin(2\eta)}{2\eta} < 1,$$

for $\theta \neq 0$ and $\eta \in (0, \pi/2)$, there are no other solutions. When $\theta = 0$, $\operatorname{Re}(m(\theta - i\eta)^2)$ is $\tan(\eta)/\eta > 1$ for all $\eta \in (0, \pi/2)$. We can thus conclude that $A \neq 0$ on $\Gamma_1$.

Similar computations yield

$$A\Big|_{\Gamma_2} = \begin{cases} 1 - m(\theta + i\eta), & \xi_0 > 0, \\ 1, & \xi_0 = 0, \end{cases}$$

which is nowhere vanishing by the same argument. Along the other arcs,

$$A\Big|_{\Gamma_3} = A\Big|_{\Gamma_4} = 1.$$

So, $A$ along $\Gamma$ is nowhere vanishing and Proposition A.1 gives the desired conclusion. □

The next result is about the Fredholm index of $\tilde{\mathcal{T}}$.

**Lemma 2.8** $\tilde{\mathcal{T}} = \tilde{t}(x, D) \colon H^j \to H^j$ *has Fredholm index two.*

***Proof*** According to Proposition A.1, the total increase of the argument of $A$ as $\Gamma$ is traversed with the counter-clockwise orientation

$$\begin{array}{ccc} (0, -1, 0, 1) & \xleftarrow{\Gamma_3} & (0, 1, 0, 1) \\ \downarrow{\Gamma_2} & & \Gamma_1 \uparrow \\ (0, -1, 0, -1) & \xrightarrow{\Gamma_4} & (0, 1, 0, -1) \end{array}$$

determines the Fredholm index of $t(x, D)$. We begin with $\Gamma_1$ from $(0, 1, 0, -1)$ to $(0, 1, 0, 1)$ and consider the total increase of the argument of $1 - m(\theta - i\eta)$ from $\theta = -\infty$ to $\theta = \infty$; see the proof of Lemma 2.7. As before, it is easier to deal with $m(\theta - i\eta)^2$. The sign of the real part of $m(\theta - i\eta)^2$ is determined by the sign of

$$\theta \sinh(2\theta) + \eta \sin(2\eta)$$

which is positive for $\theta \in \mathbb{R}$ and $\eta \in (0, \pi/2)$. This means that $m^2$ stays in the first and fourth quadrant of $\mathbb{C}$. The sign of the imaginary part equals the sign of

$$\eta \sinh(2\theta) - \theta \sin(2\eta),$$

which is a strictly increasing function in $\theta$ taking the value zero at $\theta = 0$. This means that at $\theta = 0$, $m^2$ enters the first quadrant from the fourth. By computing the value of $m^2$ as $\theta \to -\infty$, at $\theta = 0$ and as $\theta \to \infty$, we can conclude that $m^2$ along $\Gamma_1$ makes one counter-clockwise revolution about 1. Taking the square root of $m^2$ preserves the signs of the real and imaginary part. Then, multiplication with $-1$ flips the signs and addition with 1 corresponds to horizontal translation to the right by 1; see Fig. 2. Finally, we arrive at the conclusion that the increase of the argument of $A$ along $\Gamma_1$ from $(0, 1, 0, -1)$ to $(0, 1, 0, 1)$ is $2\pi$.

A similar analysis shows that an additional increase of $2\pi$ is gained along $\Gamma_2$. On $\Gamma_3$ and $\Gamma_4$, $A \equiv 1$, so there is no contribution from these arcs. In total, the increase

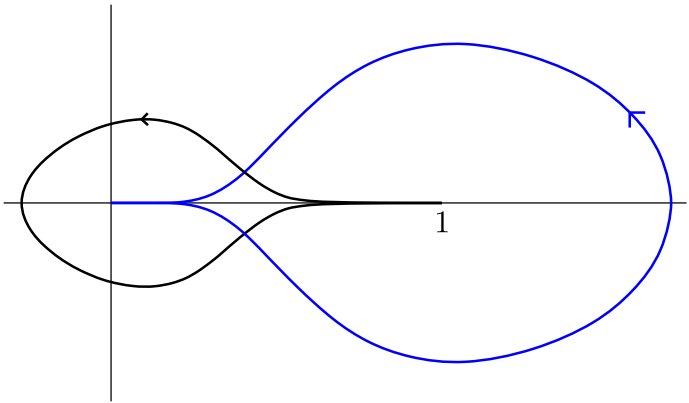

**Fig. 2** Graphs of $m^2$ (larger loop) and $1-m$ (smaller loop) along $\Gamma_1$. The increase of the argument of $1-m$ about the origin equals the increase of the argument of $m^2$ around 1, which is $2\pi$

of the argument along $\Gamma$ is $4\pi$. Proposition A.1 now gives that the Fredholm index is two. □

Conjugation with the invertible linear operator $M_{\cosh}$ preserves Fredholmness and the Fredholm index. Hence, $\mathcal{T} = M_{\cosh} \circ \tilde{\mathcal{T}} \circ M_{\cosh}^{-1} \colon H_{-\eta}^{j} \to H_{-\eta}^{j}$ is Fredholm with Fredholm index two. We have proved the first part of the main result of this section, which is the following.

**Proposition 2.9** $\mathcal{T} \colon H_{-\eta}^{j} \to H_{-\eta}^{j}$ *is Fredholm with Fredholm index two and*

$$\operatorname{Ker} \mathcal{T} = \operatorname{span}\{1, x\}.$$

*Proof* The statement concerning the Fredholm properties of $\mathcal{T}$ is already proved. Note that solving $\mathcal{T}\varphi = 0$ for $\varphi \in L_{-\eta}^{2}$ using the Fourier transform is problematic because $\mathcal{F}\varphi$ is not necessarily a tempered distribution. Thus, we consider the $L^2$-adjoint of $\mathcal{T} \colon L_{-\eta}^{2} \to L_{-\eta}^{2}$, which is $\mathcal{T} \colon L_{\eta}^{2} \to L_{\eta}^{2}$, and determine its range. The equation $\mathcal{T}\psi = g$ in $L_{\eta}^{2}$ corresponds to $(1-m)\mathcal{F}\psi = \mathcal{F}g$ on the Fourier side, where $\mathcal{F}\psi$ and $\mathcal{F}g$ are analytic functions bounded on the strip $|\operatorname{Im} z| < \eta$. In view of (3), $1 - m(\xi)$ vanishes to second order at $\xi = 0$ and is bounded away from zero if $\xi$ is. As a consequence, the range of $\mathcal{T}$ on $L_{\eta}^{2}$ consists of functions $g$ satisfying $\mathcal{F}g(0) = (\mathcal{F}g)'(0) = 0$, or equivalently $\int_{\mathbb{R}} g(x)\,dx = \int_{\mathbb{R}} xg(x)\,dx = 0$. This immediately implies that $\operatorname{Ker} \mathcal{T}$ in $L_{-\eta}^{2}$ is $\operatorname{span}\{1, x\} \subset H_{-\eta}^{j}$ and the claim is established. □

### 2.4.2 Fredholmness of $L[\varphi^*, c^*] \colon H^j \to H^j$

The application of Proposition A.1 is simpler for $L[\varphi^*, c^*] \colon H^j \to H^j$, as conjugation does not take place.

**Proposition 2.10** *Let $\varphi^*$ be a solution with wave speed $c^* > 1$, satisfying $\sup_{x \in \mathbb{R}} \varphi^*(x) < c^*/2$. Then $L[\varphi^*, c^*] \colon H^j \to H^j$ is Fredholm with Fredholm index zero.*

**Proof** Proposition A.1 is employed once again. The linear operator $L[\varphi^*, c^*]\colon \phi \mapsto c^*\phi - K * \phi - 2\varphi^*\phi$ has the symbol

$$l(x, \xi) = c^* - m(\xi) - 2\varphi^*(x) \in C^\infty(\mathbb{R} \times \mathbb{R}),$$

where smoothness of $\varphi^*$ is from Proposition 2.1(iv). The corresponding positively homogeneous function $B$ is

$$B(x_0, x, \xi_0, \xi) = l\left(\frac{x}{x_0}, \frac{\xi}{\xi_0}\right), \quad \text{where } x_0 > 0 \text{ and } \xi_0 > 0.$$

As before, we verify that $B$ does not vanish at any point along $\Gamma = \cup_{1 \le i \le 4} \Gamma_i$; see the proof of Lemma 2.7. Along $\Gamma_1$ and $\Gamma_2$,

$$\begin{aligned}
B(x_0, x, \xi_0, \xi)\Big|_{\Gamma_{1,2} \setminus \{\xi_0 = 0\}} &= \lim_{x_0 \to 0\pm} B(x_0, 1, \xi_0, \xi) \\
&= \lim_{x_0 \to 0\pm} l\left(\frac{1}{x_0}, \frac{\xi}{\sqrt{1-\xi^2}}\right) \\
&= c^* - m\left(\frac{\xi}{\sqrt{1-\xi^2}}\right),
\end{aligned}$$

as $\varphi^*$ is smooth and $\lim_{|t| \to \infty} \varphi^*(t) = 0$. Since $m \le 1$, $B$ cannot attain the value zero. Along $\Gamma_3$ and $\Gamma_4$, $B$ is $c^* - 2\varphi^*(x/\sqrt{1-x^2})$. Since $\sup_{x \in \mathbb{R}} \varphi^*(x) < c^*/2$ by assumption, $B$ cannot take the value zero. Moreover, the argument of $B$ is constant along $\Gamma$ as $B$ is real-valued and the claim is proved. □

## 3 Local bifurcation

We apply an adaptation of the nonlocal center manifold theorem in [24] to (1) in order to construct a small-amplitude solitary-solution curve emanating from $(\varphi, c) = (0, 1)$. For convenience, we work with a different bifurcation parameter $\nu := c - 1$ which will be small and positive along the local curve. In the notation of [24], Eq. (1) becomes

$$\mathcal{T}\varphi + \mathcal{N}(\varphi, \nu) = 0, \tag{10}$$

where $\mathcal{T}$ is defined in Sect. 2.4, and

$$\mathcal{N}\colon (\varphi, \nu) \mapsto \nu\varphi - \varphi^2.$$

Equation (10) will be studied for $\nu \in (0, \infty)$, in the Sobolev space of even functions $H^3_{even}$ and the weighted Sobolev spaces $H^3_{-\eta}$. This regularity choice $j = 3$ is with regard to Proposition 2.1. A solution $\varphi \in H^3_{even}$ with wave speed $c > 1$, such that $\sup_{x \in \mathbb{R}} \varphi(x) < c/2$, is smooth on $\mathbb{R}$. Moreover, $\varphi' < 0$ on $(0, \infty)$ and $\varphi$ has exponential decay.

The center-manifold reduction technique gives a reduced equation equivalent to the nonlocal Eq. (10) near the bifurcation point $(\varphi, \nu) = (0, 0)$ in $H_u^3 \times \mathbb{R}$, where $H_u^3$ is the space of functions which are uniformly local $H^3$. Since the reduced equation is an ODE, standard arguments give the existence of small-amplitude solitary-wave solutions in $H_{-\eta}^3$. Hence, by the exponential decay of solitary-wave solutions with supercritical wave speed, these are of class $H^3 \subset H_u^3$. We also prove that the bifurcation curve of non-trivial even solitary solutions is locally unique in $H_u^3 \times (0, \infty)$, and refer to it as $\mathcal{C}_{loc}$.

The global bifurcation theorem demands $L[\varphi^*, \nu^*]$ to be invertible in $H_{even}^3$ where $(\varphi^*, \nu^*) \in \mathcal{C}_{loc}$ and $\nu^* = c^* - 1$. Seeing that the Fredholm index of $L[\varphi^*, \nu^*]$ is zero, it suffices to show that the nullspace of $L[\varphi^*, \nu^*]$ is trivial. We consider Eq. (10), together with the linearized equation $L[\varphi^*, \nu^*]\phi = 0$, and formulate a center manifold theorem for this system. Exploiting the previous reduction for (10), we simplify the reduced equation for the linearized problem and are able to solve it completely in $H_{even}^3$.

### 3.1 Center manifold reductions

Two center manifold reductions are presented: one for the nonlinear problem (10) and the other for the linearized problem $L[\varphi^*, \nu^*]\phi = 0$.

For (10), we use an adaptation of the center manifold theorem in [24]. In this reference, it is assumed that the convolution kernel belongs to $W_\eta^{1,1}$, which is not the case for the Whitham kernel $K$, as $K'$ is not locally $L^1$ according to (4). Seeing that this requirement is only used for proving the Fredholm properties of the linear part $\mathcal{T}$, we replace it with requirements on the Fredholm properties on $\mathcal{T}$; see Hypothesis B.1(ii) in Appendix B.1. The rest of the proof of the center manifold theorem in [24] remains the same.

We consider (10) together with the modified equation

$$\mathcal{T}\varphi + \mathcal{N}(\chi^\delta(\varphi), \nu) = 0, \tag{11}$$

where $\chi^\delta(\varphi)$ is a nonlocal and translation invariant cutoff operator defined in Appendix B. We have $\chi^\delta(\varphi) = \varphi$ if $\|\varphi\|_{H_u^3} \leq C_0 \delta$ and $\chi^\delta(\varphi) = 0$ if $\|\varphi\|_{H_u^3}$ is sufficiently large. Since $H_u^3 \subset H_{-\eta}^3$ for all $\eta > 0$, the operator $\chi^\delta$ is a cutoff in the $H_{-\eta}^3$ norm. More details are provided in Appendix B.

We have shown that Ker $\mathcal{T}$ has dimension two in $H_{-\eta}^3$ and equals $\text{span}\{1, x\}$. Hence, elements $A + Bx \in \text{Ker }\mathcal{T}$ will often be identified with $(A, B) \in \mathbb{R}^2$. We define a projection on Ker $\mathcal{T}$,

$$\mathcal{Q} \colon \varphi \mapsto \varphi(0) + \varphi'(0)x, \quad H_{-\eta}^3 \to \text{Ker }\mathcal{T}, \tag{12}$$

which could also be considered as a mapping from $H_{-\eta}^3$ to $\mathbb{R}^2$. Finally, the shift $\varphi \mapsto \varphi(\cdot + \xi)$ will be denoted by $\tau_\xi$.

**Theorem 3.1** *For equation* (10), *there exist a neighborhood* $\mathcal{V}$ *of* $0 \in \mathbb{R}$, *a cutoff radius* $\delta$, *a weight* $\eta^* \in (0, \pi/2)$ *and a map*

$$\Psi: \mathbb{R}^2 \times \mathcal{V} \to \mathrm{Ker}\, \mathcal{Q} \subset H^3_{-\eta^*}$$

*with the center manifold*

$$\mathcal{M}^\nu_0 = \{A + Bx + \Psi(A, B, \nu) \mid A, B \in \mathbb{R}, \nu \in \mathcal{V}\} \subset H^3_{-\eta^*}$$

*as its graph. We have*

 (i) *(smoothness)* $\Psi \in \mathscr{C}^3$;
 (ii) *(tangency)* $\Psi(0, 0, 0) = 0$ *and* $\mathrm{D}_{(A,B)}\Psi(0, 0, 0) = 0$;
 (iii) *(global reduction)* $\mathcal{M}^\nu_0$ *consists precisely of* $\varphi$ *such that* $\varphi \in H^3_{-\eta^*}$ *is a solution of the modified equation* (11) *with parameter* $\nu$;
 (iv) *(local reduction) any* $\varphi$ *solving* (10) *with parameter* $\nu$ *and* $\|\varphi\|_{H^3_u} < C_0 \delta$ *is contained in* $\mathcal{M}^\nu_0$;
 (v) *(correspondence) any element* $\varphi = A + Bx + \Psi(A, B, \nu) \in \mathcal{M}^\nu_0$ *solves the local equation*

$$\varphi''(x) = f(\varphi(x), \varphi'(x), \nu), \quad \text{where } f(A, B, \nu) = \Psi''(A, B, \nu)(0), \quad (13)$$

*and conversely, any solution of this equation is an element in* $\mathcal{M}^\nu_0$. *The Taylor expansion of* $\Psi$ *gives*

$$\varphi'' = -6\varphi^2 + \frac{19}{5}(\varphi')^2 + 6\nu\varphi + \mathcal{O}\left(|(\varphi, \varphi')|(\nu^2 + |\varphi|^2 + |\varphi'|^2)\right). \quad (14)$$

 (vi) *(equivariance) besides the translations* $\tau_\xi$, *equations* (10) *and* (11) *possess a reflection symmetry* $R\varphi(x) := \varphi(-x)$, *meaning* $\mathcal{T}R\varphi = R\mathcal{T}\varphi$, $\mathcal{N}(R\varphi, \nu) = R\mathcal{N}(\varphi, \nu)$ *and* $\chi^\delta(R\varphi) = R\chi^\delta(\varphi)$. *It is hence reversible. The function* $f$ *in* (v) *commutes with all translations and anticommutes with the reflection symmetry.*

*Proof* We use Theorem B.5. Proposition 2.9 shows that Hypothesis B.1 is met. Also, the fact that $\mathcal{N}$ is a Nemytskii operator verifies Hypothesis B.3. In particular, $\mathcal{N} \in \mathscr{C}^\infty$ and $\mathcal{N}$ commutes with the translations $\tau_\xi$ for all $\nu \in (0, \infty)$. This means that in Hypothesis B.3, we can choose any regularity $k \geq 2$, possibly at the price of a smaller cutoff radius $\delta$ and weight $\eta^*$. Since a quadratic-order Taylor expansion of $\Psi$ suffices for our purposes, $k = 3$ is chosen. Statement (vi) concerning the reflection symmetry $R$ follows directly from $K$ being an even function. Hence, Theorem B.5 applies and gives items (i)–(vi).

Equation (14) in (v) is given by Theorem B.5(vii). We use $\mathcal{Q}$ defined in (12) to compute the reduced vector field. Let $\varphi \in \mathcal{M}^\nu_0$. According to Theorem B.5(viii), $\mathcal{M}^\nu_0$ is invariant under translation symmetries. Hence, $\tau_\xi \varphi$ is also an element of $\mathcal{M}^\nu_0$ for all $\xi \in \mathbb{R}$. Applying Theorem B.5(vii) on $\tau_\xi \varphi$ gives

$$\left.\frac{\mathrm{d}}{\mathrm{d}t} \mathcal{Q}(\tau_t \tau_\xi \varphi)\right|_{t=0} = \varphi'(\xi) + \varphi''(\xi)x.$$

We compute $\varphi''(\xi)$ by noting that

$$\varphi''(\xi) = \varphi''(x+\xi)\big|_{x=0} = \frac{d^2}{dx^2}(\tau_\xi\varphi)\big|_{x=0},$$

and since $\tau_\xi\varphi \in \mathcal{M}_0^\nu$,

$$\tau_\xi\varphi(x) = \varphi(\xi) + \varphi'(\xi)x + \Psi(\varphi(\xi), \varphi'(\xi), \nu)(x).$$

Hence,

$$\varphi''(\xi) = f(\varphi(\xi), \varphi'(\xi), \nu), \text{ where } f(A, B, \nu) = \Psi''(A, B, \nu)(0),$$

which is (13).

To prove Eq. (14), we compute the Taylor expansion of $\Psi$. In view of $\Psi(0, 0, 0) = 0$ and $D_{(A,B)}\Psi(0, 0, 0) = 0$, the Taylor expansion of $\Psi \colon \mathbb{R}^2 \times \mathcal{V} \to \operatorname{Ker} \mathcal{Q}$ is

$$\begin{aligned}\Psi(A, B, \nu) = &\, g(\nu)\Psi_{001} + A\nu\Psi_{101} + B\nu\Psi_{011} \\ &+ A^2\Psi_{200} + AB\Psi_{110} + B^2\Psi_{020} \\ &+ \mathcal{O}((|A| + |B|)(\nu^2 + |A|^2 + |B|^2)),\end{aligned}$$

where $\Psi_{ijk} \in \operatorname{Ker} \mathcal{Q}$ for $1 \leq i+j+k \leq 2$. Since $\mathcal{N}(0, \nu) = 0$, we have $\Psi(0, 0, \nu) = 0$ and consequently $g(\nu) = 0$ for all $\nu \in \mathbb{R}$. So, elements $\varphi$ in the center manifold $\mathcal{M}_0^\nu$ have the form

$$\begin{aligned}\varphi(x) = &\, A + Bx + \Psi(A, B, \nu)(x) \\ = &\, A + Bx + A\nu\Psi_{101}(x) + B\nu\Psi_{011}(x) \\ &+ A^2\Psi_{200}(x) + AB\Psi_{110}(x) + B^2\Psi_{020}(x) \\ &+ \mathcal{O}((|A| + |B|)(\nu^2 + |A|^2 + |B|^2)).\end{aligned}$$

Substituting $\varphi = A + Bx + \Psi(A, B, \nu)$ into (10), then identifying coefficients of orders $\mathcal{O}(A\nu)$, $\mathcal{O}(B\nu)$, $\mathcal{O}(A^2)$, $\mathcal{O}(AB)$ and $\mathcal{O}(B^2)$, we are led to the linear equations

$$\begin{aligned}\mathcal{T}\Psi_{200} &= -\mathcal{T}\Psi_{101} = 1, \\ \mathcal{T}\Psi_{110} &= -2\mathcal{T}\Psi_{011} = 2x, \\ \mathcal{T}\Psi_{020} &= x^2,\end{aligned}$$

uniquely determined by the condition $\mathcal{Q}(\Psi_{ijk}) = 0$, $i + j + k > 1$. Indeed, if there are two solutions $\Psi_{ijk}$ and $\tilde{\Psi}_{ijk}$, then $\Psi_{ijk} - \tilde{\Psi}_{ijk}$ lies in $\operatorname{Ker} \mathcal{T} \cap \operatorname{Ker} \mathcal{Q}$, and hence is zero. The linear equations are solved in Appendix C. This gives

$$\Psi(A, B, \nu)(x) = \left(-3A^2 + \frac{19}{10}B^2 + 3A\nu\right)x^2 + \left(-2AB + B\nu\right)x^3 \\ - \frac{B^2}{2}x^4 + \mathcal{O}((|A| + |B|)(\nu^2 + |A|^2 + |B|^2)).$$

Differentiating $\Psi$ twice with respect to $x$ and evaluating at $x = 0$ shows equation (14). □

When solving the linearized problem $L[\varphi^*, \nu^*]\phi = 0$, we want to take advantage of the assumption that $\varphi^* \in \mathcal{M}_0^{\nu^*}$ is a solution of (10) with parameter $\nu^*$. Hence, we consider (10) and $L[\varphi^*, \nu^*]\phi = 0$ simultaneously:

$$\mathbf{T}(\varphi, \phi) + \mathbf{N}(\varphi, \phi, \nu) = 0, \tag{15}$$

where $\mathbf{T} : (H^3_{-\eta})^2 \to (H^3_{-\eta})^2$ is an onto Fredholm operator with Fredholm index four given by

$$\mathbf{T}: (\varphi, \phi) \mapsto (\mathcal{T}\varphi, \mathcal{T}\phi),$$

and

$$\mathbf{N}: (\varphi, \phi, \nu) \mapsto (\mathcal{N}(\varphi, \nu), D_\varphi \mathcal{N}(\varphi, \nu)\phi).$$

The modified system is

$$\mathbf{T}(\varphi, \phi) + \mathbf{N}^\delta(\varphi, \phi, \nu) = 0,$$

with the nonlinearity

$$\mathbf{N}^\delta(\varphi, \phi, \nu) = (\mathcal{N}^\delta(\varphi, \nu), D_\varphi \mathcal{N}^\delta(\varphi, \nu)\phi).$$

Since we only cut off in $\varphi$, the modified linearized equation coincides with the original one, as long as $\varphi \in \mathcal{M}_0^\nu$ is sufficiently small in the $H^3_u$ topology. Hence, all solutions to the linearized equation will be captured. The downside of this scheme is that our previous adaptation of the results in [24] cannot be applied directly. We replace the contraction principle with a fiber contraction principle to achieve the following result; see Appendix B.2.

**Theorem 3.2** *For* (15), *there exist a cutoff radius* $\delta$, *a neighborhood* $\mathcal{V}$ *of* $0 \in \mathbb{R}$, *a weight* $\eta^* \in (0, \pi/2)$, *two mappings* $\Psi_1$ *and* $\Psi_2$, *where*

$$\Psi_1 \colon \mathbb{R}^2 \times \mathcal{V} \to \operatorname{Ker} \mathcal{Q} \subset H^3_{-\eta^*},$$

*with the center manifold*

$$\mathbf{M}^\nu_{0,1} := \big\{ A + Bx + \Psi_1(A, B, \nu) \,\big|\, A, B \in \mathbb{R}, \nu \in \mathcal{V} \big\}$$

*as its graph, and at each fixed element* $\varphi \in \mathbf{M}^\nu_{0,1}$ *uniquely determined by* $(A, B)$,

$$\Psi_2[A, B, \nu] : \operatorname{Ker} \mathcal{T} \to \operatorname{Ker} \mathcal{Q},$$

*with graph*

$$\mathbf{M}_{0,2}[A, B, \nu] := \{C + Dx + \Psi_2[A, B, \nu](C, D) \mid C, D \in \mathbb{R}\}.$$

*The following statements hold.*

(i) $\mathbf{M}_{0,1}^{\nu}$ *coincides with* $\mathcal{M}_0^{\nu}$ *in Theorem* 3.1 *and all statements in this theorem hold for* $\mathbf{M}_{0,1}^{\nu}$;

(ii) $\Psi_2[A, B, \nu] = D_{(A,B)}\Psi_1(A, B, \nu)$, *so* $\Psi_2[A, B, \nu]$ *is a bounded linear operator from* Ker $\mathcal{T}$ *to* Ker $\mathcal{Q}$. *Also,* $\Psi_2$ *is* $\mathscr{C}^{k-1}$ *in* $(A, B, \nu)$. *Suppose that* $\varphi^* \in \mathbf{M}_{0,1}^{\nu^*}$ *is sufficiently small in the* $H_u^3$ *norm, so that* $\varphi^*$ *is a solution of* (10) *with parameter* $\nu^*$, *uniquely determined by* $(A^*, B^*)$. *Then,*

(iii) $\mathbf{M}_{0,2}[A^*, B^*, \nu^*]$ *consists precisely of solutions* $\phi \in H_{-\eta^*}^3$ *of the linearized equation* $\mathcal{T}\phi + D_{\varphi}\mathcal{N}(\varphi^*, \nu^*)\phi = 0$;

(iv) *every* $\phi \in \mathbf{M}_{0,2}[A^*, B^*, \nu^*]$ *is a solution of*

$$\phi''(x) = g(\varphi^*(x), (\varphi^*)'(x), \phi(x), \phi'(x), \nu^*),$$

*where*

$$g(A^*, B^*, C, D, \nu^*) = D_A f(A^*, B^*, \nu^*)C + D_B f(A^*, B^*, \nu^*)D,$$

*and* $f(A, B, \nu) = \Psi_1''(A, B, \nu)(0)$. *Conversely, any solution of the above second-order ODE is an element in* $\mathbf{M}_{0,2}[A^*, B^*, \nu^*]$. *The Taylor expansion of g gives*

$$\phi'' = -12\varphi^*\phi + \frac{38}{5}(\varphi^*)'\phi' + 6\nu^*\phi \tag{16}$$
$$+ \mathcal{O}((|\varphi^*\phi| + |(\varphi^*)'\phi'|)((\nu^*)^2 + |\varphi^*| + |(\varphi^*)'|)).$$

**Proof** Theorem B.6 applies and gives (i)–(iv). The cutoff radii $\delta$ given by Theorem 3.1 and Theorem 3.2 are not necessarily the same, but the smallest one can be chosen to have (i). Arguing along the same lines as the proof of Theorem 3.1 gives equation (iv) and differentiating the Taylor expansion of $f$ in (14) gives Eq. (16). □

### 3.2 Local bifurcation curve

The center manifold theorem states that a solution $\varphi$ sufficiently small in the $H_u^3$ norm solves the reduced ODE (14), which is local in nature and allows spatial dynamics tools. Let $\varphi = P$, $\varphi' = Q$ and regard the spatial variable $x$ as "time" $t$. Equation (14) defines the following system of ODEs

$$\begin{cases} \dfrac{dP}{dt} = Q \\ \dfrac{dQ}{dt} = -6P^2 + \frac{19}{5}Q^2 + 6\nu P + \mathcal{O}((|P| + |Q|)(\nu^2 + |P|^2 + |Q|^2)), \end{cases} \tag{17}$$

which is reversible by Theorem 3.1(vi). We aim to rescale (17) into a KdV-equation when $\nu = 0$, that is

$$\begin{cases} \dfrac{d\tilde{P}}{dT} = \tilde{Q}(T) \\ \dfrac{d\tilde{Q}}{dT} = \tilde{P}(T) - \dfrac{3}{2}\tilde{P}(T)^2. \end{cases}$$

Hence, we set

$$T = \alpha t, \quad P(t) = \beta \tilde{P}(T), \quad Q(t) = \gamma \tilde{Q}(T).$$

Differentiating, substituting into (17) and identifying coefficients yield

$$\frac{\gamma}{\beta\alpha} = 1, \quad 6\nu\frac{\beta}{\alpha\gamma} = 1, \quad \frac{6\beta^2}{\alpha\gamma} = \frac{3}{2},$$

which are satisfied by

$$\alpha = \sqrt{6\nu}, \quad \beta = \frac{3}{2}\nu, \quad \gamma = \sqrt{\frac{3^3\nu^3}{2}}.$$

The resulting rescaled system is

$$\begin{cases} \dfrac{d\tilde{P}}{dT} = \tilde{Q}(T) \\ \dfrac{d\tilde{Q}}{dT} = \tilde{P}(T) - \dfrac{3}{2}\tilde{P}(T)^2 + \dfrac{57}{10}\nu\tilde{Q}(T)^2 + \mathcal{O}\left(\nu(|\tilde{P}| + \nu^{1/2}|\tilde{Q}|)(1 + |\tilde{P}|^2 + \nu|\tilde{Q}|^2)\right). \end{cases} \quad (18)$$

For $\nu = 0$, (18) is the KdV-equation with the explicitly known pair of solutions

$$\tilde{P}(T) = \operatorname{sech}^2\left(\frac{T}{2}\right), \quad \tilde{Q}(T) = -\operatorname{sech}^2\left(\frac{T}{2}\right)\tanh\left(\frac{T}{2}\right),$$

which corresponds to a symmetric and homoclinic orbit. For $\nu > 0$, the symmetric homoclinic orbit persists by the same argument as in [26], p. 955. Undoing the rescaling and switching back to $P = \varphi$ as well as $Q = \varphi'$ give

$$\begin{aligned} \varphi(t) &= \frac{3}{2}\nu \operatorname{sech}^2\left(\frac{\sqrt{6\nu}t}{2}\right) + \mathcal{O}(\nu^2) \\ \varphi'(t) &= -\frac{3^{3/2}\nu^{3/2}}{2^{1/2}}\operatorname{sech}^2\left(\frac{\sqrt{6\nu}t}{2}\right)\tanh\left(\frac{\sqrt{6\nu}t}{2}\right) + \mathcal{O}(\nu^{5/2}), \end{aligned} \quad (19)$$

for $\nu > 0$. The supercritical solitary-wave solution $\varphi$ is exponentially decaying, so both $\varphi$ and $\varphi'$ belong to $H^3 \subset H^3_u$. Also, they depend continuously on the parameter

$\nu$. We denote this solution as $\varphi^*_{\nu^*}$ with parameter $\nu^*$ and define

$$\mathcal{C}_{loc} = \{(\varphi^*_{\nu^*}, \nu^*) \mid 0 < \nu^* < \nu'\},$$

for some $\nu' > 0$. The main result of this section is reached.

**Theorem 3.3** *There exists a neighborhood of $(\varphi, \nu) = (0, 0)$ in $H^3_u \times (0, \infty)$, for which $\mathcal{C}_{loc}$ is the unique $\nu$-dependent family of non-trivial even small-amplitude solitary solutions to (10) emanating from $(0, 0)$. We refer to $\mathcal{C}_{loc}$ as the local bifurcation curve.*

*Proof* The function $\varphi^*_{\nu^*}$ belongs to $\mathcal{M}^{\nu^*}_0$ by the one-to-one correspondence between (17) and $\mathcal{M}^{\nu^*}_0$ in Theorem 3.1(v). From (19) combined with the fact that $\varphi^*_{\nu^*}$ is exponentially decaying, we have $\|\varphi^*_{\nu^*}\|_{H^1_u} \lesssim \nu^*$. Since the reduced vector field $f$ in (13) is superlinear in $\varphi^*_{\nu^*}(x)$ and $(\varphi^*_{\nu^*})'(x)$ by Theorem 3.1(ii), the bound by $\nu^*$ in (19) is carried over to $(\varphi^*_{\nu^*})''$. Differentiating (13) gives

$$(\varphi^*_{\nu^*})^{(3)} = D_1 f(\varphi^*_{\nu^*}, (\varphi^*_{\nu^*})', \nu^*) \cdot (\varphi^*_{\nu^*})' + D_2 f(\varphi^*_{\nu^*}, (\varphi^*_{\nu^*})', \nu^*) \cdot (\varphi^*_{\nu^*})'',$$

where $D_1 f$ and $D_2 f$ are bounded in view of Theorem 3.1(i). Hence, $(\varphi^*_{\nu^*})^{(3)}$ is also bounded by $\nu^*$. We obtain the improvement $\|\varphi^*_{\nu^*}\|_{H^3_u} \lesssim \nu^*$ and then by choosing $\nu^*$ sufficiently small, $\varphi^*_{\nu^*}$ is indeed a solution of (10) according to Theorem 3.1(iv). The existence of $\mathcal{C}_{loc}$ in $H^3_{even}$ is now established since $\varphi^*_{\nu^*} \in H^3$ is an even function.

Our argument for the uniqueness of $\mathcal{C}_{loc}$ is similar to the one in [11], Lemma 5.10. Suppose that $(\varphi, \nu)$ is a non-trivial even solitary wave solution which is small enough in $H^3_u \times (0, \infty)$ that $\varphi$ lies on the center manifold $\mathcal{M}^\nu_0$. Then $(P, Q) = (\varphi, \varphi')$ is a reversible homoclinic solution of the ODE (17), whose phase portrait is qualitatively the same as in Fig. 3. The homoclinic orbit in the right half plane corresponds to the case $\varphi = \varphi_\nu$ and hence $(\varphi, \nu) \in \mathcal{C}_{loc}$. Any other solution must therefore approach the origin along the portions of its stable and unstable manifolds lying in the left half plane. But this would force $P(t) = \varphi(t) < 0$ for sufficiently large $|t|$, contradicting (i) in Proposition 2.1. □

*Remark 3.4* Since the amplitude of $\varphi^*_{\nu^*}$ is $\mathcal{O}(\nu^*)$ as $\nu^* \to 0$, we can find a $\nu' > 0$, such that

$$\sup_{x \in \mathbb{R}} \varphi^*_{\nu^*}(x) < \frac{1 + \nu^*}{2}, \quad \text{for all } \nu^* \in (0, \nu').$$

From Proposition 2.1, the solutions $\varphi^*_{\nu^*}$ are everywhere smooth and strictly decreasing on $(0, \infty)$.

### 3.3 Invertibility of $L[\varphi^*, \nu^*]$

Let $\varphi^* := \varphi^*_{\nu^*} \in \mathcal{C}_{loc}$ and the corresponding parameter $\nu^*$ be fixed. The linear operator $L[\varphi^*, \nu^*]$ is the linearization of the left-hand side of (10) with respect to the $\varphi$-component. Note that $L[\varphi^*, \nu^*]: H^3_{even} \to H^3_{even}$ is Fredholm with Fredholm index zero. Hence, we only need to show that the nullspace of $L[\varphi^*, \nu^*]$ is trivial. The invertibility of $L[\varphi^*, \nu^*]$ has already been shown in [31]. In this section, we showcase an alternative approach exploiting Theorem 3.2 and are able to make quantitative statements for elements in the nullspace of $L[\varphi^*, \nu^*]$ in $H^3_{-\eta^*}$, where $\eta^*$ is as in Theorem 3.2. This approach is inspired by Lemmas 4.14 and 4.15 in [32].

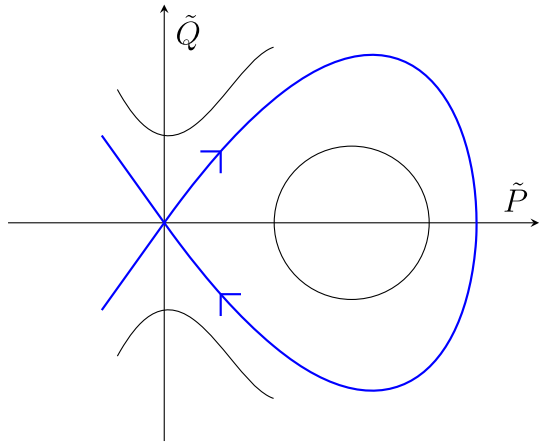

**Fig. 3** Phase portrait for (18) when $\nu = 0$. The homoclinic orbit persists for $\nu > 0$

**Proposition 3.5** *The nullspace* Ker $L$ *of* $L[\varphi^*, \nu^*] \colon H^3_{-\eta^*} \to H^3_{-\eta^*}$ *is two-dimensional, spanned by the exponentially decaying function* $(\varphi^*)'$ *and an exponentially growing function. Seeing that* $(\varphi^*)'$ *is odd,* Ker $L$ *restricted on* $H^3_{\text{even}}$ *is trivial and* $L[\varphi^*, \nu^*]$ *is thus invertible in* $H^3_{\text{even}}$.

***Proof*** We use Theorem 3.2(iv), that is, elements $\phi \in \text{Ker } L \subset H^3_{-\eta^*}$ have a one-to-one correspondence to the solutions of (16). Letting $\phi = U$, $\phi' = V$ and regarding $x$ as a time variable $t$, we cast (16) into a system

$$\begin{cases} \frac{dU}{dt} = V \\ \frac{dV}{dt} = 6\nu^* U - 12\varphi^* U + \frac{38}{5}(\varphi^*)' V + \mathcal{O}((|\varphi^* U| + |(\varphi^*)' V|)((\nu^*)^2 + |\varphi^*| + |(\varphi^*)'|)). \end{cases} \quad (20)$$

This can be considered as a perturbation problem of the form

$$\frac{d\mathbf{u}}{dt} = (M + R(t))\mathbf{u},$$

where $\mathbf{u} \colon \mathbb{R}^2 \to \mathbb{R}^2$, $M \in \mathbb{R}^{2 \times 2}$ is a matrix with constant coefficients and $R \colon \mathbb{R} \to \mathbb{R}^2$ is an integrable remainder term. In this case,

$$M = \begin{pmatrix} 0 & 1 \\ 6\nu^* & 0 \end{pmatrix},$$

so the eigenvalues for $M$ are $\sqrt{6\nu^*}$ and $-\sqrt{6\nu^*}$. Moreover, the exponential decay of $\varphi^*$ and $(\varphi^*)'$ guarantees the integrability condition; see the discussion after (19).

Applying for example Problem 29, Chapter 3 in [13], as well as switching back to $\phi$ and $\phi'$, the statements concerning Ker $L$ are immediate. In particular, Ker $L$ is spanned by a function $\phi_1$ behaving as $\exp(t\sqrt{6\nu^*})$ and $\phi_2$ behaving as $\exp(-t\sqrt{6\nu^*})$ as $t \to \infty$. It is a straightforward calculation to show that one exponentially decaying function in Ker $L$ is $(\varphi^*)'$, which is an odd function. Since even functions in $H^3$ cannot be written as linear combinations of an odd exponentially decaying function and an exponentially growing function, Ker $L$ is trivial in $H^3_{even}$. □

## 4 Global bifurcation

We use a global bifurcation theorem from [11] in a slightly modified form because the open set in our case is not a product set; see Appendix D. For (1), we take

$$\mathcal{X} = \mathcal{Y} = H^3_{even}, \quad F(\varphi, \nu) = \mathcal{T}\varphi + \mathcal{N}(\varphi, \nu),$$

and

$$\mathcal{U} = \left\{ (\varphi, \nu) \in \mathcal{X} \times (0, \infty) \,\middle|\, \sup_{x \in \mathbb{R}} \varphi(x) < \frac{1+\nu}{2} \right\}.$$

Since $H^3 \subset BUC^2$, the supremum norm is controlled by the $H^3$ norm and $\mathcal{U}$ is thus an open set in $H^3_{even} \times \mathbb{R}$. We aim to use Theorem D.1. Proposition 2.10 verifies Hypothesis (A) in this theorem, while Sect. 3.2 and Proposition 3.5 together verify Hypothesis (B). Here, the local curve $\mathcal{C}_{loc}$ bifurcates from $(0, 0) \in \partial \mathcal{U}$. We have thus the following global bifurcation theorem for (1) in $H^3_{even}$ and $\mathcal{U}$.

**Theorem 4.1** *The local bifurcation curve $\mathcal{C}_{loc}$ in Section 3.2 is contained in a curve of solutions $\mathcal{C}$, which is parametrized as*

$$\mathcal{C} = \{(\varphi_s, \nu_s) \mid 0 < s < \infty\} \subset \mathcal{U} \cap F^{-1}(0)$$

*for some continuous map $(0, \infty) \ni s \mapsto (\varphi_s, \nu_s)$. We have*

(a) *One of the following alternatives holds:*

  (i) *(blowup) as $s \to \infty$,*

  $$M(s) := \|\varphi_s\|_{H^3} + \nu_s + \frac{1}{\text{dist}((\varphi_s, \nu_s), \partial \mathcal{U})} \to \infty;$$

  (ii) *(loss of compactness) there exists a sequence $s_n \to \infty$ as $n \to \infty$ such that $\sup_n M(s_n) < \infty$ but $(\varphi_{s_n})_n$ has no subsequence convergent in $\mathcal{X}$.*

(b) *Near each point $(\varphi_{s_0}, \nu_{s_0}) \in \mathcal{C}$, we can reparametrize $\mathcal{C}$ so that $s \mapsto (\varphi_s, \nu_s)$ is real analytic.*

(c) *$(\varphi_s, \nu_s) \notin \mathcal{C}_{loc}$ for $s$ sufficiently large.*

In this section, we use the integral identity in Proposition 2.2 to exclude the loss of compactness scenario. An alternative route is to employ the Hamiltonian structures for nonlocal problems in [7]. Even though the Whitham kernel does not fit into this framework, a direct differentiation confirms that equation 43 in [7] indeed gives a Hamiltonian for the Whitham equation. We also study how $M(s)$ blows up as $s \to \infty$.

### 4.1 Preservation of nodal structure

We begin by showing that the nodal structure is preserved along the global bifurcation curve.

**Theorem 4.2** *If $(\varphi, \nu) \in \mathcal{C} \subset \mathcal{U}$, then $\varphi$ is smooth on $\mathbb{R}$ and strictly decreasing on the interval $(0, \infty)$.*

*Proof* The property $\sup_{x \in \mathbb{R}} \varphi(x) < (1 + \nu)/2$ for $(\varphi, \nu) \in \mathcal{U}$ implies smoothness on $\mathbb{R}$ by Proposition 2.1(iv). Because in addition $\varphi \in H^3_{even}$ is a solitary-wave solution with parameter $\nu > 0$, we have that $\varphi$ is non-increasing by Proposition 2.1(iii). In order to apply Proposition 2.1(v) to conclude that $\varphi$ is strictly decreasing on $(0, \infty)$ we only need to establish that $\varphi$ is non-constant.

The only constant solutions are 0 and $\nu > 0$ and the latter is excluded by the fact that $\varphi \in H^3$ is a solitary-wave solution. To show that $\varphi(x) \not\equiv 0$, note that the linearization $L[0, \nu] \colon \phi \mapsto (1 + \nu)\phi - K * \phi$ on $H^3_{even}$ is Fredholm of index zero for all $\nu \in \mathcal{I}$ by Proposition 2.10 and Ker $L[0, \nu]$ is trivial. So, $L[0, \nu]$ is invertible. The implicit function theorem applies and prevents $\mathcal{C}$ from intersecting the trivial solution line. Hence, this alternative cannot occur. □

### 4.2 Compactness of the global curve $\mathcal{C}$

The following result rules out alternative (ii) in Theorem 4.1(a).

**Theorem 4.3** *Every sequence $(\varphi_n, \nu_n)_{n=1}^\infty := (\varphi_{s_n}, \nu_{s_n})_{n=1}^\infty \subset \mathcal{C}$ satisfying*

$$\sup_n M(s_n) < \infty$$

*has a convergent subsequence in $H^3_{even} \times (0, \infty)$.*

*Proof* Proposition 2.3 gives a locally uniform limit $\varphi$ for a subsequence of functions $\varphi_n$. The idea is to use Proposition 2.5 to show that $\varphi_n$ converges to $\varphi$ in $H^3$. We observe that the assumption $\sup_n M(s_n) < \infty$ implies

$$\sup_n \|\varphi_n\|_{H^3} < \infty \text{ and } \inf_n \nu_n > 0,$$

where $\nu_n = c_n - 1$. Also, according to Corollary 2.4, $\varphi$ inherits non-negativity, continuity, evenness, boundedness and monotonicity from $\varphi_n$. More precisely, we have shown in Theorem 4.2 that $\varphi_n$ is strictly decreasing on $(0, \infty)$. So, $\varphi$ is at least non-increasing on $(0, \infty)$.

First, we verify that $\varphi_n \to \varphi$ uniformly. Because the sequence of functions $\varphi_n$ is uniformly bounded in $H^3$, it has a weak limit which coincides with the locally uniform limit $\varphi$. So $\varphi \in H^3$. Since $\varphi$ is in addition monotone on the real half-lines, we have $\lim_{|x|\to\infty} \varphi(x) = 0$. Corollary 2.4(ii) now confirms the desired uniform convergence of $\varphi_n$ to $\varphi$.

Next, for the $L^2$ convergence, we use the integral identity in Proposition 2.2. For each $n$, $\varphi_n \in BUC^2$ and $\lim_{|x|\to\infty} \varphi_n(x) = 0$. Hence,

$$\lim_{R\to\infty} \int_{-R}^{R} \varphi_n(\varphi_n - \nu_n) \, dx = 0.$$

Also, since $\nu_n > 0$, the solitary-wave solution $\varphi_n$ has exponential decay according to Proposition 2.1(ii) and we are allowed to write

$$\int_{\mathbb{R}} \varphi_n^2 \, dx = \nu_n \int_{\mathbb{R}} \varphi_n \, dx.$$

Since $\sup_n \|\varphi_n\|_{H^3} < \infty$, the $L^2$ integral on the left-hand side is uniformly bounded in $n$. Because $\inf_n \nu_n > 0$, the $L^1$ integral on the right-hand side is uniformly bounded as well. Taking into account that $\lim_{|x|\to\infty} \varphi_n(x) = 0$ uniformly in $n$, we obtain

$$\int_{\mathbb{R}} \varphi_n^2 \, dx = \int_{|x|<R_\epsilon} \varphi_n^2 \, dx + \int_{|x|>R_\epsilon} \varphi_n^2 \, dx$$

$$\leq \int_{|x|<R_\epsilon} \varphi_n^2 \, dx + \epsilon \cdot \sup_n \|\varphi_n\|_{L^1}$$

As $n \to \infty$, we have $\int_{|x|<R_\epsilon} \varphi_n^2 \, dx \to \int_{|x|<R_\epsilon} \varphi^2 \, dx$. Letting $\epsilon \to 0$ confirms that $\varphi_n \to \varphi$ in $L^2$.

Finally, we observe that $\sup_n M(s_n) < \infty$ also implies

$$\sup_{x\in\mathbb{R}} \varphi_n(x) < (1+\nu_n)/2 \quad \text{and} \quad \sup_{x\in\mathbb{R}} \varphi(x) < (1+\nu)/2.$$

All prerequisites of Proposition 2.5 are now checked and we have $\varphi_n \to \varphi$ in $H^3$. □

### 4.3 Analysis of the blowup

Having excluded the loss of compactness alternative, we examine the blowup alternative

$$\lim_{s\to\infty} \left( \|\varphi_s\|_{H^3} + \nu_s + \frac{1}{\text{dist}((\varphi_s, \nu_s), \partial\mathcal{U})} \right) = \infty,$$

where $(\varphi_s, \nu_s) \in \mathcal{C}$. In this case, for any sequence $s_n \to \infty$ we can extract a subsequence (also denoted $\{s_n\}$) for which at least one of the following four possibilities holds:

(P1) $\|\varphi_{s_n}\|_{H^3} \to \infty$,
(P2) $\nu_{s_n} \to \infty$,
(P3) $\nu_{s_n} \to 0$,
(P4) $(1 + \nu_{s_n})/2 - \sup_{x \in \mathbb{R}} \varphi_{s_n}(x) \to 0$,

where (P3) and (P4) belong to the case when $\text{dist}((\varphi_{s_n}, \nu_{s_n}), \partial \mathcal{U}) \to 0$.

**Theorem 4.4** *The alternatives (P2) and (P3) cannot occur.*

**Proof** Alternative (P2) cannot occur since the definition of $\mathcal{U}$ and Proposition 2.1(viii) imply that $\nu_{s_n} \leq 1$.

To exclude alternative (P3), we assume $\nu_{s_n} \to 0$ as $n \to \infty$. Any locally uniform limit $(\varphi, 0)$ solves (1). Moreover, $\varphi$ is bounded, continuous, and monotone; see Proposition 2.3. Proposition 2.2 gives

$$\inf_x \varphi(x) < 0 < \sup_x \varphi(x) \text{ or } \varphi \equiv 0.$$

Because $\varphi$ is non-negative, we must have $\varphi \equiv 0$. In particular, $\lim_{|x| \to \infty} \varphi(x) = 0$. In virtue of Corollary 2.4(ii), $\varphi_{s_n} \to 0$ uniformly and now according to Remark 2.6, $\varphi_{s_n} \to 0$ in $C^k$ for any $k$, which implies that $\varphi_{s_n} \to \varphi$ in $H^3_{-\eta}$ for all $\eta > 0$ and that $(\varphi_{s_n}, \nu_{s_n})$ reenters any small neighborhood of $(0, 0)$ in $H^3_u \times (0, \infty)$. This cannot happen in light of Theorem 4.1(c) and the uniqueness of $\mathcal{C}_{loc}$ given by Theorem 3.3. □

Next, we show a useful characterization for when the $H^3$ norm stays bounded.

**Lemma 4.5** *By possibly taking a subsequence of $(\varphi_{s_n}, \nu_{s_n})_{n=1}^\infty$, we have $\nu_{s_n} \to \nu > 0$ and $\varphi_{s_n} \to \varphi$ locally uniformly as $n \to \infty$. Then, $\inf_n \nu_{s_n} > 0$. Moreover,*

$$\sup_n \|\varphi_{s_n}\|_{H^3} < \infty \text{ if and only if } \sup_{x \in \mathbb{R}} \varphi(x) < \frac{1 + \nu}{2} \text{ and } \lim_{|x| \to \infty} \varphi(x) = 0.$$

**Proof** The existence of such a subsequence is given by Proposition 2.3. Since Theorem 4.4 has excluded (P3), we must have $\nu_{s_n} \to \nu > 0$, which also implies that $\inf_n \nu_{s_n} > 0$. We focus on proving the last statement. The proof of Theorem 4.3 already gives that $\sup_n \|\varphi_{s_n}\|_{H^3} < \infty$ implies $\varphi \in H^3$, thus $\sup_{x \in \mathbb{R}} \varphi(x) < (1+\nu)/2$ by Proposition 2.1(v) and $\lim_{|x| \to \infty} \varphi(x) = 0$. Conversely, assume on the contrary that there is a subsequence of functions $\varphi_{s_n}$ such that $\|\varphi_{s_n}\|_{H^3} \to \infty$ as $n \to \infty$, yet its locally uniform limit $\varphi$ satisfies $\sup_{x \in \mathbb{R}} \varphi(x) < (1+\nu)/2$ and $\lim_{|x| \to \infty} \varphi(x) = 0$. Then, $\varphi$ is smooth by (iv) in Proposition 2.1. Also, Corollary 2.4(ii) gives $\varphi_{s_n} \to \varphi$ uniformly, so $\varphi_{s_n}(x) \to 0$ as $|x| \to \infty$ uniformly in $n$. Similar to the proof of Theorem 4.3, we get

$$\nu_{s_n} \int_{\mathbb{R}} \varphi_{s_n} \, dx = \int_{\mathbb{R}} \varphi_{s_n}^2 \, dx \leq \int_{|x| < R_\epsilon} \varphi_{s_n}^2 \, dx + \epsilon \int_{\mathbb{R}} \varphi_{s_n} \, dx,$$

where $\epsilon$ and $R_\epsilon$ are independent of $n$. Rearranging gives

$$\int_{\mathbb{R}} \varphi_{s_n} \, dx \leq \frac{1}{v_{s_n} - \epsilon} \int_{|x| < R_\epsilon} \varphi_{s_n}^2 \, dx \leq \frac{2R_\epsilon}{v_{s_n} - \epsilon} \cdot \sup_{x \in \mathbb{R}} \varphi_{s_n}^2(x) \leq \frac{2R_\epsilon}{v_{s_n} - \epsilon},$$

where $\sup_{x \in \mathbb{R}} \varphi_{s_n}^2(x) < 1$ because $v_{s_n} \in (0, 1]$. Recall that $\inf_n v_{s_n} > 0$. Choosing $\epsilon = \inf_n v_{s_n}/2$, this shows $\sup_n \|\varphi_{s_n}\|_{L^1} < \infty$. It follows that the sequence of functions $\varphi_{s_n}$ is uniformly bounded in $L^2$ and arguing as in the proof of Proposition 2.5 gives the uniform boundedness in $H^3$, which is a contradiction to the assumption. □

We can now establish the following equivalence.

**Theorem 4.6** *(P4) and (P1) are equivalent.*

***Proof*** Let $(\varphi_{s_n}, v_{s_n})_{n=1}^\infty$ be a sequence satisfying (P4). By possibly taking a subsequence, Proposition 2.3 gives that $\varphi_{s_n} \to \varphi$ locally uniformly and $v_{s_n} \to v$, where $v > 0$ as we have excluded (P3). Since each $\varphi_{s_n}$ is even and strictly decreasing, (P4) is the same as

$$\lim_{n \to \infty} \left| \frac{1 + v_{s_n}}{2} - \varphi_{s_n}(0) \right| = 0,$$

which is equivalent to

$$\varphi(0) = \frac{1 + v}{2}.$$

By Lemma 4.5, this implies (P1). For the other implication, let $(\varphi_{s_n}, v_{s_n})_{n=1}^\infty$ be a sequence satisfying (P1). We also have that $\varphi_{s_n} \to \varphi$ locally uniformly and $v_{s_n} \to v > 0$. Once again by Proposition 2.3 and Corollary 2.4, $\varphi$ solves (1) with parameter $v > 0$ and is continuous, bounded, even, and $\varphi$ is non-increasing on $(0, \infty)$. Then, the limit $\lim_{|x| \to \infty} \varphi(x)$ exists. According to Corollary 2.4(i), this can take the value

$$\text{(a)} \lim_{|x| \to \infty} \varphi(x) = 0 \quad \text{or} \quad \text{(b)} \lim_{|x| \to \infty} \varphi(x) = v > 0.$$

In addition, Proposition 2.2 says

$$\text{(A)} \inf_{x \in \mathbb{R}} \varphi(x) < v < \sup_{x \in \mathbb{R}} \varphi(x), \quad \text{(B)} \; \varphi \equiv 0 \quad \text{or} \quad \text{(C)} \; \varphi \equiv v > 0.$$

The combinations (aC) and (bB) are quickly excluded. If $\varphi \equiv 0$, then $\sup_{x \in \mathbb{R}} \varphi(x) < (1+v)/2$ and Lemma 4.5 gives that $\sup_n \|\varphi_{s_n}\|_{H^3} < \infty$, which contradicts the blowup alternative. This rules out (aB). The fact that $\varphi$ is non-increasing on $(0, \infty)$ rules out (bA). Assume (bC), which is just (C). Then, we arrive at a contradiction as follows. Consider the sequence of translated solutions $\tau_{x_n} \varphi_{s_n} = \varphi_{s_n}(\cdot + x_n)$, where each $x_n$ is chosen in such a way that

$$\varphi_{s_n}(x_n) = \tilde{v}, \quad \text{where} \quad 0 < \tilde{v} < \inf_n \frac{v_{s_n}}{2}.$$

Such a number $\tilde{\nu}$ exists because $\nu_{s_n} > 0$ cannot limit to 0. Moreover, $\lim_{n\to\infty} x_n = \infty$, or we cannot have $\varphi \equiv \nu > 0$ for every $x$ while each $\varphi_n$ has exponential decay. Corollary 2.4(iii) applied to $(\tau_{x_n}\varphi_{s_n})_n$ gives a bounded, continuous and non-increasing locally uniform limit $\tilde{\varphi}$. The function $\tilde{\varphi}$ is a solution to (10) with parameter $\nu$. Its limits $\tilde{\varphi}(x)$ as $x \to \pm\infty$ are guaranteed to exist and these can take the value zero or $\nu > 0$. By construction,

$$\tilde{\varphi}(0) = \lim_n \varphi_{s_n}(0 + x_n) = \tilde{\nu} \in (0, \nu),$$

which implies

$$\lim_{x\to-\infty} \tilde{\varphi}(x) = \nu \text{ and } \lim_{x\to\infty} \tilde{\varphi}(x) = 0.$$

On the other hand, this also shows that $\tilde{\varphi}$ is not a constant function, implying

$$0 = \inf_{x\in\mathbb{R}} \tilde{\varphi}(x) < \nu < \sup_{x\in\mathbb{R}} \tilde{\varphi}(x) = \nu,$$

which is a contradiction. We conclude that (C) cannot occur. Hence, we must have

$$\lim_{|x|\to\infty} \varphi(x) = 0 \text{ and } \inf_{x\in\mathbb{R}} \varphi(x) < \nu < \sup_{x\in\mathbb{R}} \varphi(x),$$

where the first condition is the same as (P4). Applying Lemma 4.5 gives the desired implication. □

**Remark 4.7** This in fact implies that $\mathcal{C}$ satisfies $\|\varphi_s\|_{H^3} \to \infty$ and $(1 + \nu_s)/2 - \sup_{x\in\mathbb{R}} \varphi_s(x) \to 0$ as $s \to \infty$, without the need for considering subsequences.

Finally, since (P1) and (P4) are equivalent, (P2) and (P3) cannot happen and the blowup alternative must take place, we must have a sequence $(\varphi_{s_n}, \nu_{s_n})_{n=1}^\infty \subset \mathcal{C}$ such that $\lim_{n\to\infty} \|\varphi_{s_n}\|_{H^3} = \infty$ and $\lim_{n\to\infty} \nu_{s_n} = \nu > 0$. By taking the limit of a subsequence, an extreme solitary-wave solution $\varphi$ attaining the highest possible amplitude $(1 + \nu)/2$ is found; see Fig. 1.

**Theorem 4.8** *There exists a sequence of elements $(\varphi_{s_n}, \nu_{s_n}) \in \mathcal{C}$, such that*

$$\lim_{n\to\infty} \|\varphi_{s_n}\|_{H^3} = \infty \quad \text{and} \quad \lim_{n\to\infty} \nu_{s_n} = \nu > 0.$$

*The sequence of solutions $\varphi_{s_n}$ has a locally uniform limit $\varphi$. We have*

(i) *$\varphi$ is continuous, bounded, even and non-increasing on the positive real half-line;*
(ii) *$\varphi$ is a non-trivial solitary-wave solution to (10) with parameter $\nu > 0$;*
(iii) *$\varphi(0) = (1 + \nu)/2$ and more precisely*

$$C_1|x|^{\frac{1}{2}} \leq \frac{1+\nu}{2} - \varphi(x) \leq C_2|x|^{\frac{1}{2}}, \tag{21}$$

*near the origin and for some constants $0 < C_1 < C_2$;*

(iv) $\varphi$ *is smooth everywhere except at* $x = 0$;
(v) $\varphi$ *has exponential decay.*

**Proof** Existence has already been shown. According to Proposition 2.3, there is such a locally uniform limit $\varphi$. Statement (i) is immediate from Corollary 2.4(ii). Statement (ii) follows from the proof of Theorem 4.6, as we have shown that

$$\lim_{|x|\to\infty} \varphi(x) = 0 \text{ and } \inf_{x\in\mathbb{R}} \varphi(x) < \nu < \sup_{x\in\mathbb{R}} \varphi(x).$$

Statement (iii) is Theorem 4.6 and the estimate (21) for $\varphi$ near the origin is Proposition 2.1(vi). Together with $\varphi$ being non-increasing on $(0,\infty)$, we have $\varphi(x) < (1+\nu)/2$ if $x \neq 0$. Proposition 2.1(iv) applies and gives statement (iv). Statement (v) follows from Proposition 2.1(ii). □

**Acknowledgements** This research is supported by the Swedish Research Council, grant no. 2016-04999. Part of it was carried out during the workshop *Nonlinear water waves—an interdisciplinary interface* in 2017 at Erwin Schrödinger International Institute for Mathematics and Physics. The hospitality and support of the institute is gratefully acknowledged. We also thank Grégory Faye and Arnd Scheel for the swift and helpful correspondence during the finishing phase of this project. Finally, we thank the referee for helpful comments.

**Funding.** Open access funding provided by Lund University.

# Appendix

## A Fredholmness of pseudodifferential operators

Let $x^* = (x_0, x) \in \mathbb{R}^2$ and

$$X^* = \{x^* \in \mathbb{R}^2 \mid x_0 \geq 0, x^* \neq 0\}.$$

Similarly, let $\xi^* = (\xi_0, \xi) \in \mathbb{R}^2$ and

$$E^* = \{\xi^* \in \mathbb{R}^2 \mid \xi_0 \geq 0, \xi^* \neq 0\}.$$

We define

$$\mathbb{S} = \{(x_0, x, \xi_0, \xi) \in \mathbb{R}^2 \times \mathbb{R}^2 \mid x_0^2 + |x|^2 = \xi_0^2 + |\xi|^2 = 1, x_0 > 0, \xi_0 > 0\},$$

and denote the relative closure of $\mathbb{S}$ in $X^* \times E^*$ by $\overline{\mathbb{S}}$ and the boundary of $\overline{\mathbb{S}}$ by $\Gamma$.

Let $\mathcal{A}$ be the class of functions $A(x^*, \xi^*) \in C^\infty(X^* \times E^*)$ such that $A$ is positively homogeneous of degree 0 in $x^*$ and $\xi^*$, that is,

$$A(\lambda x^*, \xi^*) = A(x^*, \lambda \xi^*) = A(x^*, \xi^*), \quad \lambda > 0.$$

Clearly, each $A \in \mathcal{A}$ is uniquely determined by its values on $\mathbb{S}$. Conversely, each function $\tilde{A} \in C^\infty(\overline{\mathbb{S}})$ can be uniquely homogeneously extended to $X^* \times E^*$. So, $\mathcal{A} \cong C^\infty(\overline{\mathbb{S}})$. By $S^0_{\mathcal{A}}$, we denote the set of symbols $p_A(x, \xi)$ which are given by

$$p_A(x, \xi) = A(1, x, 1, \xi),$$

for some $A \in \mathcal{A}$. For $p_A \in S^0_{\mathcal{A}}$, we have the following result, which combines Theorems 4.1 and 4.2 in [27].

**Proposition A.1** *Let $j \in \mathbb{R}$. If $p_A(x, \xi) \in S^0_{\mathcal{A}}$ and $A(x^*, \xi^*) \neq 0$ on $\Gamma$, then*

$$p_A(x, D) \colon H^j \to H^j$$

*is Fredholm with Fredholm index*

$$\operatorname{ind} p_A(x, D) = \frac{1}{2\pi} \left( \arg A(x^*, \cdot^*) \Big|_0 \right),$$

*where* $\arg A(x^*, \cdot^*)|_0$ *is the increase in the argument of $A(x^*, \xi^*)$ around $\Gamma$ as $\Gamma$ is traversed with the counter-clockwise orientation.*

We comment that a version of Proposition A.1 is available for matrix-valued symbols $p_A(x, \xi)$; see [27, Sect. 4].

## B Center manifold theorems

### B.1 A nonlocal version of the center manifold theorem

Nonlocal nonlinear parameter-dependent problems of the form

$$\mathcal{T}v + \mathcal{N}(v, \mu) = 0, \tag{22}$$

where

$$\mathcal{T}v = v + \mathcal{K} * v,$$

are considered in the weighted Sobolev spaces $H^j_{-\eta}(\mathbb{R}^n)$ for some $\eta > 0$, positive integers $n$ and $j = 1$. The following presentation focuses on dimension $n = 1$ and arbitrary regularity $j \geq 1$. $\mathcal{T}$ will be referred to as the linear part and $\mathcal{N}$ as the nonlinear part of (22).

To obtain an appropriate modification of Eq. (22), we introduce the space of uniformly local $H^j$ functions,

$$H_{\mathrm{u}}^j = \left\{ v \in H_{\mathrm{loc}}^j \,\middle|\, \|v\|_{H_{\mathrm{u}}^j} < \infty \right\} \text{ where } \|v\|_{H_{\mathrm{u}}^j} = \sup_{y \in \mathbb{R}} \|v(\,\cdot\, + y)\|_{H^j([0,1])}.$$

Note that the embeddings $H^j \subset H_{\mathrm{u}}^j \subset H_{-\eta}^j$ are continuous for all $\eta > 0$. Next, let $\underline{\chi}: \mathbb{R} \to \mathbb{R}$ be a smooth cutoff function satisfying

$$\underline{\chi}(x) = \begin{cases} 1, & |x| < 1 \\ 0, & |x| > 2, \end{cases} \quad \text{and } \sup_{x \in \mathbb{R}} |\underline{\chi}'(x)| \leq 2.$$

In addition, let $\theta : \mathbb{R} \to \mathbb{R}$ be a smooth and even function with

$$\sum_{j \in \mathbb{Z}} \theta(x - j) = 1, \quad \operatorname{supp} \theta \subset [-1, 1], \quad \theta(x) \geq 0, \quad \theta\left(\left[0, \tfrac{1}{2}\right]\right) \subset \left[\tfrac{1}{2}, 1\right],$$

so that

$$\int_{y \in \mathbb{R}} \theta(x - y)\, \mathrm{d}y = \int_{y=0}^{1} \sum_{j \in \mathbb{Z}} \theta(x - y - j)\, \mathrm{d}y = 1.$$

We define

$$\chi : v \mapsto \int_{y \in \mathbb{R}} \underline{\chi}(\|\theta(\,\cdot\, - y)v\|_{H^j}) \theta(x - y) v(x)\, \mathrm{d}y,$$

and then

$$\chi^\delta : v \mapsto \delta \cdot \chi\left(\frac{v}{\delta}\right), \quad \delta > 0.$$

The cutoff operator $\chi \colon H_{-\eta}^j \to H_{\mathrm{u}}^j$ is well-defined, Lipschitz continuous and $\|\chi(v)\|_{H_{\mathrm{u}}^j} \leq C_0$ for all $v \in H_{-\eta}^j$. Lemma 2 in [25] proves these claims for $j = 1$ and an asymmetric function $\theta$ with $\operatorname{supp} \theta \subset [-1/4, 5/4]$. Generalizing to higher regularity $j \geq 1$ and verifying the results in [25] with our choice of $\theta$ are straightforward. The scaled cutoff operator $\chi^\delta \colon H_{-\eta}^j \to H_{\mathrm{u}}^j$ naturally inherits these properties, in particular $\|\chi^\delta(v)\| \leq C_0 \delta$ for all $v \in H_{-\eta}^j$. The modification of Eq. (22) is

$$\mathcal{T} v + \mathcal{N}^\delta(v, \mu) = 0, \tag{23}$$

where

$$\mathcal{N}^\delta : (v, \mu) \mapsto \mathcal{N}(\chi^\delta(v), \mu).$$

For $\xi \in \mathbb{R}$, $\tau_\xi v = v(\,\cdot\, + \xi)$ denotes the shift by $\xi$. Regarded as an operator on $H^j_{-\eta}$ or $H^j$, $\tau_\xi$ is bounded with $\|\tau_\xi\|_{H^j_{-\eta} \to H^j_{-\eta}} \lesssim \exp(\eta\xi)$ and $\|\tau_\xi\|_{H^j \to H^j} = 1$.

Also, we let $\mathcal{Q}\colon H^j_{-\eta} \to H^j_{-\eta}$ be a bounded projection on the nullspace Ker $\mathcal{T}$ of $\mathcal{T}$ with a continuous extension to $H^{j-1}_{-\eta}$, such that $\mathcal{Q}$ commutes with the inclusion map from $H^j_{-\eta}$ to $H^j_{-\eta'}$, for all $0 < \eta' < \eta$.

**Hypothesis B.1** (The linear part $\mathcal{T}$)
  (i) There exists $\eta_0 > 0$ such that $\mathcal{K} \in L^1_{\eta_0}$.
  (ii) The operator
$$\mathcal{T}\colon v \mapsto v + \mathcal{K} * v, \quad H^j_{-\eta} \to H^j_{-\eta}$$
is Fredholm for $\eta \in (0, \eta_0)$, its nullspace Ker $\mathcal{T}$ is finite-dimensional and
$$\operatorname{ind} \mathcal{T} = \dim \operatorname{Ker} \mathcal{T},$$
where ind $\mathcal{T}$ is the Fredholm index of $\mathcal{T}$. In other words, $\mathcal{T}$ is onto.

**Remark B.2** Hypothesis B.1(i) gives boundedness of $\mathcal{T}$ on $H^j_{-\eta}$ for $\eta \in (0, \eta_0)$; see Sect. 2.4. The original assumption $\mathcal{K} \in W^{1,1}_\eta$ in [24] is only used to guarantee item (ii). Since the Whitham kernel $K \notin W^{1,1}_\eta$, we instead require (ii) directly.

**Hypothesis B.3** (The nonlinear part $\mathcal{N}$) There exist $k \geq 2$, a neighborhood $\mathcal{U}$ of $0 \in H^j_{-\eta}$ and $\mathcal{V}$ of $0 \in \mathbb{R}$, such that for all sufficiently small $\delta > 0$, we have

  (i) $\mathcal{N}^\delta\colon H^j_{-\eta} \times \mathcal{V} \to H^j_{-\eta}$ is $\mathscr{C}^k$. In addition, for all non-negative pairs $(\zeta, \eta)$ satisfying $0 < k\zeta < \eta < \eta_0$, $D^l_v \mathcal{N}^\delta\colon (H^j_{-\zeta})^l \to H^j_{-\eta}$ is bounded for all $0 < l\zeta \leq \eta < \eta_0$ and $0 \leq l \leq k$, and Lipschitz continuous in $v$ for $1 \leq l \leq k-1$ uniformly in $\mu \in \mathcal{V}$;
  (ii) $\mathcal{N}^\delta(\tau_\xi v, \mu) = \tau_\xi \mathcal{N}^\delta(v, \mu)$ for all $\mu \in \mathcal{V}$ and $\xi \in \mathbb{R}$;
  (iii) $\mathcal{N}^\delta(0,0) = 0$, $D_v \mathcal{N}^\delta(0,0) = 0$ and as $\delta \to 0$,
$$\delta_1(\delta) := \operatorname{Lip}_{H^j_{-\eta} \times \mathcal{V}} \mathcal{N}^\delta = \mathcal{O}(\delta + |\mu|).$$

In dimension $n \geq 1$, a symmetry is a triple $(\rho, \tau_\xi, \kappa) \in \mathbf{O}(n) \times (\mathbb{R} \times \mathbf{O}(1))$, where the orthogonal linear transformation $\rho \in \mathbf{O}(n)$ acts on $v(x) \in \mathbb{R}^n$ while $\tau_\xi$ and $\kappa$ act on the variable $x \in \mathbb{R}$. In particular, a symmetry $(\rho, \tau_\xi, \kappa)$ is called equivariant if $(\rho, \tau_\xi, \kappa) \in \mathbf{O}(n) \times (\mathbb{R} \times \{\mathrm{Id}\})$, and reversible otherwise. Lemma 2 in [25] and the fact that $\theta$ is even give the invariance of $\chi^\delta$ under the whole group $\mathbf{O}(n) \times (\mathbb{R} \times \mathbf{O}(1))$ for $n = 1$, that is, $\chi^\delta(\gamma v) = \gamma \chi^\delta(v)$ for all $\gamma \in \mathbf{O}(1) \times (\mathbb{R} \times \mathbf{O}(1))$.

**Hypothesis B.4** *(Symmetries)* There exists a symmetry group $S \subset \mathbf{O}(1) \times (\mathbb{R} \times \mathbf{O}(1))$ which contains all translations on the real line and which commutes with the linear part $\mathcal{T}$ as well as the nonlinear part $\mathcal{N}$, that is,
$$\mathcal{T}(\gamma v) = \gamma(\mathcal{T} v) \text{ and } \mathcal{N}(\gamma v, \mu) = \gamma \mathcal{N}(v, \mu), \text{ for all } \gamma \in S.$$

We state the equivariant parameter-dependent center manifold theorem. Even though the hypothesis on $\mathcal{K}$ is changed, the proof is the same.

**Theorem B.5** *Assume Hypotheses* B.1, B.3 *and* B.4 *are met for equation* (22). *Then, by possibly shrinking the neighborhood $\mathcal{V}$ of $0 \in \mathbb{R}$, there exist a cutoff radius $\delta > 0$, a weight $\eta^* \in (0, \eta_0)$ and a map*

$$\Psi \colon \operatorname{Ker} \mathcal{T} \times \mathcal{V} \subset H^j_{-\eta^*} \times \mathbb{R} \to \operatorname{Ker} \mathcal{Q} \subset H^j_{-\eta^*}$$

*with the center manifold*

$$\mathcal{M}_0^\mu := \left\{ v_0 + \Psi(v_0, \mu) \,\big|\, v_0 \in \operatorname{Ker} \mathcal{T}, \mu \in \mathcal{V} \right\} \subset H^j_{-\eta^*},$$

*as its graph. The following statements hold:*

(i) *(smoothness)* $\Psi \in \mathscr{C}^k$, *where k is as in Hypothesis* B.3;
(ii) *(tangency)* $\Psi(0, 0) = 0$ *and* $\mathrm{D}_{v_0} \Psi(0, 0) = 0$;
(iii) *(global reduction)* $\mathcal{M}_0^\mu$ *consists precisely of functions v such that* $v \in H^j_{-\eta^*}$ *is a solution of the modified equation* (23) *with parameter $\mu$;*
(iv) *(local reduction) any function* $v \in H^j_u$ *solving* (22) *with* $\|v\|_{H^j_u} < C_0 \delta$ *is contained in* $\mathcal{M}_0^\mu$;
(v) *(translation invariance) the shift $\tau_\xi$ by $\xi \in \mathbb{R}$ acting on $\mathcal{M}_0^\mu$ induces a $\mu$-dependent flow*

$$\Phi_\xi \colon \operatorname{Ker} \mathcal{T} \to \operatorname{Ker} \mathcal{T}$$

*through* $\Phi_\xi = \mathcal{Q} \circ \tau_\xi \circ (\operatorname{Id} + \Psi)$;
(vi) *(reduced vector field) the reduced flow $\Phi_\xi$ is of class $\mathscr{C}^k$ in $v_0, \mu, \xi$ and is generated by a reduced parameter-dependent vector field $f$ of class $\mathscr{C}^{k-1}$ on the finite-dimensional* $\operatorname{Ker} \mathcal{T}$;
(vii) *(correspondence) any element $v = v_0 + \Psi(v_0, \mu)$ of $\mathcal{M}_0^\mu$ corresponds one-to-one to a solution of*

$$\frac{\mathrm{d} v_0}{\mathrm{d} x} = f(v_0, \mu) := \frac{\mathrm{d}}{\mathrm{d} x} \mathcal{Q}(\tau_x v) \Big|_{x=0}.$$

(viii) *(equivariance)* $\operatorname{Ker} \mathcal{T}$ *is invariant under $S$ and $\mathcal{Q}$ can be chosen to commute with all $\gamma \in S$. Consequently, $\Psi$ commutes with $\gamma \in S$ and $\mathcal{M}_0^\mu$ is invariant under $S$. Finally, the reduced vector field $f$ in item (vi) commutes with all equivariant symmetries and anticommutes with the reversible ones in $S$.*

### B.2 Center manifold theorem for an augmented problem

In this section, we prove a center manifold theorem for a system consisting of the nonlocal Eq. (22) and its linearization at $(v, \mu)$, that is,

$$\mathcal{T} w + \mathrm{D}_v \mathcal{N}(v, \mu) w = 0. \tag{24}$$

More precisely, we study

$$\mathbf{T}(v, w) + \mathbf{N}(v, w, \mu) = 0, \tag{25}$$

where the linear part is

$$\mathbf{T}\colon (v, w) \mapsto (\mathcal{T}v, \mathcal{T}w) = (v + \mathcal{K} * v, w + \mathcal{K} * w),$$

and the nonlinear part is

$$\mathbf{N}\colon (v, w, \mu) \mapsto (\mathcal{N}(v, \mu), \mathrm{D}_v \mathcal{N}(v, \mu)w).$$

We consider the modified system

$$\mathbf{T}(v, w) + \mathbf{N}^\delta(v, w, \mu) = 0, \tag{26}$$

with

$$\mathbf{N}^\delta(v, w, \mu) = (\mathcal{N}^\delta(v, \mu), \mathrm{D}_v \mathcal{N}^\delta(v, \mu)w),$$

where $\mathcal{N}^\delta(v, \mu)$ is defined in the previous section. Observe that we only cut off in $v$, which allows capturing all solutions of the linearized Eq. (24). This requires yet another adaptation of the center manifold theorem, where the usual contraction principle is replaced with a fiber contraction principle.

**Theorem B.6** *Assume Hypotheses* B.1, B.3, B.4 *and let $\mathcal{Q}$ be the same projection as in Theorem* B.5. *For* (25), *there exist a cutoff radius $\delta > 0$, a weight $\eta^* \in (0, \eta_0)$, a neighborhood $\mathcal{V}$ of $0 \in \mathbb{R}$ and two mappings*

$$\Psi_1 \colon \operatorname{Ker} \mathcal{T} \times \mathcal{V} \to \operatorname{Ker} \mathcal{Q} \subset H^j_{-\eta^*},$$

*with the center manifold*

$$\mathbf{M}^\mu_{0,1} := \left\{ v_0 + \Psi_1(v_0, \mu) \,\big|\, v_0 \in \operatorname{Ker} \mathcal{T}, \mu \in \mathcal{V} \right\}$$

*as its graph, and at each fixed element $v = v_0 + \Psi_1(v_0, \mu) \in \mathbf{M}^\mu_{0,1}$,*

$$\Psi_2[v_0, \mu] \colon \operatorname{Ker} \mathcal{T} \to \operatorname{Ker} \mathcal{Q},$$

*with graph*

$$\mathbf{M}_{0,2}[v_0, \mu] := \left\{ w_0 + \Psi_2[v_0, \mu](w_0) \,\big|\, w_0 \in \operatorname{Ker} \mathcal{T} \right\}.$$

*The following statements hold.*

(i) $\mathbf{M}_{0,2}[v_0, \mu]$ consists precisely of the solutions to the modified linearized equation $\mathcal{T}w + D_v \mathcal{N}^\delta(v, \mu)w = 0$ at $v = v_0 + \Psi_1(v_0, \mu) \in \mathbf{M}_{0,1}^\mu$ and the problem

$$\mathcal{T}w + D_v \mathcal{N}^\delta(v, \mu)w = 0 \text{ with } \mathcal{Q}w = w_0$$

has a unique solution for each $w_0 \in \text{Ker } \mathcal{T}$.

(ii) $\mathbf{M}_{0,1}^\mu$ coincides with $\mathcal{M}_0^\mu$ in Theorem B.5 and all statements in this theorem hold for $\mathbf{M}_{0,1}^\mu$. In particular, if $v = v_0 + \Psi_1(v_0, \mu) \in \mathbf{M}_{0,1}^\mu$ is sufficiently small in the $H_u^j$ norm, so that $v$ solves the original equation (22) with parameter $\mu$, then $w \in \mathbf{M}_{0,2}[v_0, \mu]$ solves the linearized equation (24) at $v$.

(iii) We have $\Psi_2[v_0, \mu] = D_{v_0} \Psi_1(v_0, \mu)$. Consequently, $\Psi_2[v_0, \mu]$ is a bounded linear map, $\Psi_2[0, 0] = 0$ and $\Psi_2$ is $\mathscr{C}^{k-1}$ in $(v_0, \mu)$.

(iv) The shift $\tau_\xi$ acting on $\mathbf{M}_{0,2}[\,\cdot\,, \mu]$ induces a flow

$$\Phi_{2,\xi} : \text{Ker } \mathcal{T} \to \text{Ker } \mathcal{T}$$

through $\Phi_{2,\xi} := \mathcal{Q} \circ \tau_\xi \circ (\text{Id} + \Psi_2[\tau_\xi v_0, \mu])$.

(v) The reduced flow $\Phi_{2,\xi}$ is of class $\mathscr{C}^{k-1}$ in $v_0, \mu, \xi$ and boundedly linear in $w_0$. It is generated by a reduced parameter-dependent vector field $g$ of class $\mathscr{C}^{k-2}$ in $(v_0, \mu)$ and boundedly linear in $w_0$.

(vi) Any element $w = w_0 + \Psi_2[v_0, \mu]w_0$ of $\mathbf{M}_{0,2}[v_0, \mu]$ corresponds one-to-one to a solution of

$$\frac{dw_0}{dx} = g(v_0, w_0, \mu) := \frac{d}{dx} \mathcal{Q}\big(\tau_x w_0 + \Psi_2[\tau_x v_0, \mu](\tau_x w_0)\big)\Big|_{x=0}.$$

We have $g(v_0, w_0, \mu) = D_{v_0} f(v_0, \mu) w_0$.

We will use the following fiber contraction theorem in the proof; see Sect. 1.11.8 in [12].

**Proposition B.7** *Let $\mathcal{X}$ and $\mathcal{Y}$ be complete metric spaces. Consider a continuous map $\Lambda \colon \mathcal{X} \times \mathcal{Y} \to \mathcal{X} \times \mathcal{Y}$ of the form*

$$\Lambda(x, y) = (\lambda_1(x), \lambda_2(x, y)),$$

*where $\lambda_1 \colon \mathcal{X} \to \mathcal{X}$ and $\lambda_2 \colon \mathcal{X} \times \mathcal{Y} \to \mathcal{Y}$. If $\lambda_1$ is a contraction in $\mathcal{X}$, and $y \mapsto \lambda_2(x, y)$ is a contraction in $\mathcal{Y}$ for every fixed $x \in \mathcal{X}$, then $\Lambda$ has a unique fixed point $(x_0, y_0) \in \mathcal{X} \times \mathcal{Y}$.*

*Proof of Theorem B.6* We present the necessary changes in the proof of the center manifold theorem in [24] without the parameter $\mu$. The transition to the parameter-dependent version is the same as in [24].

In view of Hypothesis B.1(ii), the operator $\mathbf{T} = (\mathcal{T}, \mathcal{T})$ is Fredholm and its Fredholm index is twice the index of $\mathcal{T}$. We define the Fredholm-bordered operator

$$\tilde{\mathbf{T}} \colon (H_{-\eta}^j)^2 \to (H_{-\eta}^j \times \text{Ker } \mathcal{T})^2$$

with

$$\tilde{\mathbf{T}}\colon (v, w) \mapsto (\tilde{\mathcal{T}}v, \tilde{\mathcal{T}}w) := (\mathcal{T}v, \mathcal{Q}v, \mathcal{T}w, \mathcal{Q}w).$$

By Lemma 3.2 in [24], $\tilde{\mathbf{T}}$ is invertible for any $\eta \in (0, \eta_0)$ with

$$\tilde{\mathbf{T}}^{-1} = (\tilde{\mathcal{T}}^{-1}, \tilde{\mathcal{T}}^{-1}), \text{ and } \|\tilde{\mathbf{T}}^{-1}\|_{H^j_{-\eta} \to H^j_{-\eta}} \leq C(\eta),$$

where $C$ is a continuous function of $\eta$. We define the bordered nonlinearity

$$\tilde{\mathbf{N}}^\delta \colon (v, w, v_0, w_0) \mapsto (\mathcal{N}^\delta(v), -v_0, D_v \mathcal{N}^\delta(v)w, -w_0).$$

Due to Hypothesis B.3, $\tilde{\mathbf{N}}^\delta$ is continuous in $(v, w) \in (H^j_{-\eta})^2$. The bordered equation becomes

$$\tilde{\mathbf{T}}(v, w) + \tilde{\mathbf{N}}^\delta(v, w, v_0, w_0) = 0.$$

Applying $\tilde{\mathbf{T}}^{-1}$ on both sides, then moving the nonlinear term to the right-hand side, we obtain a fixed point equation

$$\begin{aligned}(v, w) &= -\tilde{\mathbf{T}}^{-1}\left(\tilde{\mathbf{N}}^\delta(v, w, v_0, w_0)\right) \\ &=: \mathcal{S}^\delta(v, w, v_0, w_0).\end{aligned}$$

Let $v_0$ and $w_0$ be fixed. We apply Proposition B.7 with $\Lambda = \mathcal{S}^\delta$ and $\mathcal{X} = \mathcal{Y} = H^j_{-\eta}$. Continuity of $\mathcal{S}^\delta$ in $(v, w)$ is clear in view of Hypothesis B.3(i) and Theorem B.5 already shows that the first component of $\mathcal{S}^\delta$ is a contraction mapping for appropriately chosen cutoff radii $\delta > 0$ and weights $\eta \in (0, \eta_0)$. Now, let $v$ be fixed in the second component of $\mathcal{S}^\delta$. Hypothesis B.3(iii) together with $\mathcal{N}^\delta \colon H^j_{-\eta} \to H^j_{-\eta}$ being $\mathscr{C}^k$ for $k \geq 2$ imply that

$$\sup_{v \in H^j_{-\eta}} \|D_v \mathcal{N}^\delta(v)\|_{H^j_{-\eta} \to H^j_{-\eta}} = \mathcal{O}(\delta), \text{ as } \delta \to 0.$$

This gives

$$\|D_v \mathcal{N}^\delta(v)w_1 - D_v \mathcal{N}^\delta(v)w_2\|_{H^j_{-\eta}} \lesssim \delta \|w_1 - w_2\|_{H^j_{-\eta}}.$$

Choosing $\delta$ sufficiently small, the second component is a contraction mapping on $H^j_{-\eta}$ for each fixed $v$. By Theorem B.7, there exists a unique fixed point $(v, w) := \tilde{\Psi}(v_0, w_0)$ for each prescribed $v_0$ and $w_0$:

$$(v, w) = \tilde{\Psi}(v_0, w_0) := (v_0 + \Psi_1(v_0), w_0 + \Psi_2[v_0](w_0)),$$

where $\mathcal{Q}v = v_0$, $\mathcal{Q}w = w_0$ and $\Psi_i : \text{Ker}\,\mathcal{T} \to \text{Ker}\,\mathcal{Q}$ for $i = 1, 2$. Item (i) is thus established.

Item (ii) follows from the uniqueness of $\mathcal{M}_0$ and the fact that $\chi^\delta(v) = v$ if the $H_u^j$ norm of $v$ is sufficiently small. Since elements of $\mathbf{M}_{0,1}$ are precisely solutions of the modified equation $\mathcal{T}v + \mathcal{N}^\delta(v) = 0$, we insert $v = \mathcal{Q}v + \Psi_1(\mathcal{Q}v)$ into the equation and differentiate with respect to $v$. This gives that

$$w = \mathcal{Q}w + D\Psi_1(\mathcal{Q}v)\mathcal{Q}w = w_0 + D_{v_0}\Psi_1(v_0)w_0$$

is a solution of the modified linearized equation $\mathcal{T}w + D_v\mathcal{N}^\delta(v)w = 0$. Since $D_{v_0}\Psi_1(v_0)w_0 \in \text{Ker}\,\mathcal{Q}$, we get $\mathcal{Q}(w_0 + D_{v_0}\Psi(v_0)w_0) = w_0$. By uniqueness of a solution $w$ for each $\mathcal{Q}w = w_0$, we obtain $\Psi_2[v_0] = D_{v_0}\Psi_1(v_0)$. The remaining claims in (iii) are straightforward in view of (ii) and Theorem B.5. It follows from Hypothesis B.4 that $\tau_\xi \circ \mathcal{N}^\delta = \mathcal{N}^\delta \circ \tau_\xi$, which after a differentiation gives

$$\tau_\xi\left(D_v\mathcal{N}^\delta(v)w\right) = D_v\mathcal{N}^\delta(\tau_\xi v)\tau_\xi w.$$

Then, reasoning as in [24] validates statements (iv)–(vi), except for the last claim $g(v_0, w_0) = D_{v_0}f(v_0)w_0$. This is shown by plugging $\tau_x w = (\text{Id} + D_{v_0}\Psi_1(\tau_x v_0))\tau_x w_0$ into the reduced vector field $g$ in (vi), and then identifying the result with $D_{v_0}f = D_{v_0}(\mathcal{Q} \circ \tau_x \circ (\text{Id} + \Psi_1))$. □

## C Computation of the coefficients $\Psi_{ijk}$

We wish to solve

$$\mathcal{T}\Psi_{200} = -\mathcal{T}\Psi_{101} = 1, \tag{27}$$
$$\mathcal{T}\Psi_{110} = -2\mathcal{T}\Psi_{011} = 2x, \tag{28}$$
$$\mathcal{T}\Psi_{020} = x^2, \tag{29}$$

subjected to the condition $\mathcal{Q}(\Psi_{ijk}) = 0$ for all $i + j + k \geq 1$. This condition is imposed for unique solvability of $\Psi_{ijk}$; see the proof of Theorem 3.1.

Using the fact that multiplication by $x^n$ corresponds to $n$-times differentiation on the Fourier side, we have

$$\int_\mathbb{R} K(x)x^n\,dx = \begin{cases} 0 & \text{if } n \text{ is odd} \\ (-1)^{n/2}m^{(n)}(0) & \text{if } n \text{ is even.} \end{cases}$$

Now, we can compute the convolution of $K$ and monomials $x^n$, where $n \in \mathbb{N}$. For instance, $K * 1 = 1$ and

$$\int_\mathbb{R} K(y)(x-y)\,dy = x\int_\mathbb{R} K(y)\,dy - \int_\mathbb{R} yK(y)\,dy = x \cdot (K * 1) - 0 = x,$$

where $\int_{\mathbb{R}} y K(y) \, dy = 0$ because the integrand is odd. Utilizing the binomial theorem to expand $(x - y)^n$ together with symmetries of the integrands, we arrive at

$$x^2 - \int_{\mathbb{R}} K(y)(x - y)^2 \, dy = m''(0), \tag{30}$$

$$x^3 - \int_{\mathbb{R}} K(y)(x - y)^3 \, dy = 3m''(0)x, \tag{31}$$

$$x^4 - \int_{\mathbb{R}} K(y)(x - y)^4 \, dy = 6m''(0)x^2 - m^{(4)}(0). \tag{32}$$

To solve (27), we are motivated by (30) and make the Ansatz $\Psi_{200} = \alpha x^2 - \mathcal{Q}(\alpha x^2)$, where subtraction by $\mathcal{Q}(\alpha x^2)$ is to make sure that $\mathcal{Q}(\Psi_{200}) = 0$. Since $\mathcal{Q}(\alpha x^2) = 0$ for all $\alpha \in \mathbb{R}$, it can be removed. Plugging the Ansatz into (27) yields

$$\mathcal{T}\Psi_{200} = \alpha m''(0) = 1.$$

Since $m''(0) = -1/3$, $\alpha = -3$ necessarily. In conclusion, $\Psi_{200} = -3x^2$. The linear Eqs. (28) and (29) are solved in a similar way. We summarize the results below.

$$\Psi_{101} = 3x^2, \quad \Psi_{200} = -3x^2,$$
$$\Psi_{110} = -2x^3, \quad \Psi_{020} = -\frac{1}{2}x^4 + \frac{19}{10}x^2,$$
$$\Psi_{011} = x^3.$$

## D A global bifurcation theorem

Let $\mathcal{X}, \mathcal{Y}$ be Banach spaces and $\mathcal{U} \subset \mathcal{X} \times \mathbb{R}$ an open set. Consider the abstract operator equation

$$F(\varphi, \nu) = 0,$$

where $F : \mathcal{U} \to \mathcal{Y}$ is an analytic mapping. The following is a version of Theorem 6.1 in [11] which has been slightly modified to better fit the situation in the present paper. The proof remains the same.

**Theorem D.1** *Assume*

(A) *for all $(\varphi^*, \nu^*) \in \mathcal{U} \cap F^{-1}(0)$, the Fréchet derivative $D_\varphi F(\varphi^*, \nu^*)$ is Fredholm of index zero;*
(B) *there exists a local curve of solutions $\mathcal{C}_{loc}$ with a continuous parametrization $(0, \nu') \ni \nu^* \mapsto (\varphi^*_{\nu^*}, \nu^*)$, so that*

$$\mathcal{C}_{loc} = \{(\varphi^*_{\nu^*}, \nu^*) \mid 0 < \nu^* < \nu'\} \subset \mathcal{U} \cap F^{-1}(0)$$

*and*

$$\lim_{\nu^* \to 0^+} (\varphi^*_{\nu^*}, \nu^*) \in \partial \mathcal{U},$$

*as well as*

$$D_\varphi F(\varphi^*_{\nu^*}, \nu^*) \colon \mathcal{X} \to \mathcal{Y} \text{ is invertible for all } (\varphi^*_{\nu^*}, \nu^*) \in \mathcal{C}_{loc}.$$

*Then, $\mathcal{C}_{loc}$ is contained in a curve of solutions $\mathcal{C}$, which is parametrized as*

$$\mathcal{C} = \{(\varphi_s, \nu_s) \,|\, 0 < s < \infty\} \subset \mathcal{U} \cap F^{-1}(0)$$

*for some continuous map $(0, \infty) \ni s \mapsto (\varphi_s, \nu_s)$. The global curve $\mathcal{C}$ has the following properties.*

(a) *One of the following alternatives holds:*
  (i) *(blowup) as $s \to \infty$,*
  
  $$M(s) := \|\varphi_s\|_\mathcal{X} + |\nu_s| + \frac{1}{\text{dist}((\varphi_s, \nu_s), \partial \mathcal{U})} \to \infty;$$
  
  (ii) *(loss of compactness) there exists a sequence $s_n \to \infty$ as $n \to \infty$ such that $\sup_n M(s_n) < \infty$ but $(\varphi_{s_n})_n$ has no subsequence convergent in $\mathcal{X}$.*

(b) *Near each point $(\varphi_{s_0}, \nu_{s_0}) \in \mathcal{C}$, we can reparametrize $\mathcal{C}$ so that $s \mapsto (\varphi_s, \nu_s)$ is real analytic.*

(c) $(\varphi_s, \nu_s) \notin \mathcal{C}_{loc}$ *for $s$ sufficiently large.*